\documentclass[11pt]{article}

\newif\ifpdf
\ifx\pdfoutput\undefined
\pdffalse
\else
\pdfoutput=1
\pdftrue
\fi

\ifpdf
\usepackage[pdftex]{graphicx}
\else
\usepackage{graphicx}

\fi
\usepackage{textcomp}
\usepackage{amsmath}
\usepackage{amsfonts}
\usepackage{epsfig}

\renewcommand{\baselinestretch} {1.3}

\makeatletter
\def\singlespace{\def\baselinestretch{1}\@normalsize}

\title{{\sc A Functional Wavelet-Kernel Approach for Continuous-time Prediction}}

\author{  {Anestis ~{\sc ANTONIADIS}},\\
          Laboratoire IMAG-LMC,\\
          University Joseph Fourier,\\
          BP 53, 38041 Grenoble Cedex 9,\\
          FRANCE.\\
          Email: \texttt{Anestis.Antoniadis@imag.fr}
          \\
          \\
          {Efstathios ~{\sc PAPARODITIS}}\\
          Department of Mathematics and Statistics, \\
          University of Cyprus,\\
          P.O. Box 20537,\\
          CY 1678 Nicosia,\\
          CYPRUS.\\
          Email: \texttt{stathisp@ucy.ac.cy}
          \\
          \\
          and
          \\
          \\
          {Theofanis ~{\sc SAPATINAS}}\\
          Department of Mathematics and Statistics, \\
          University of Cyprus,\\
          P.O. Box 20537,\\
          CY 1678 Nicosia,\\
          CYPRUS.\\
          Email: \texttt{t.sapatinas@ucy.ac.cy}
}

\date{}

\newcommand{\figh}[3]{\begin{figure}[htbp]       %
\begin{center}

\ifpdf
\includegraphics[width=#2]{#1}
\else
\includegraphics[width=#2]{#1.eps}
\fi

\end{center}

        \small                                   %
        \begin{singlespace}                      %
        \caption{#3}                             %
        \label{#1}                               %
        \end{singlespace}                        %
        \end{figure}}                            %

\newcommand{\reals}{\ensuremath{{\mathbb R}}}
\newcommand{\RR}{\reals}

\newcommand{\NN}{\ensuremath{{\mathbb N}}}
\newcommand{\EE}{\ensuremath{{\mathbb E}}}
\newcommand{\PP}{\ensuremath{{\mathbb P}}}

\newcommand{\KK}{\ensuremath{{\mathbb K}}}

\newtheorem{theorem}{Theorem}[section]

\newtheorem{remark}{Remark}[section]

\begin{document}
\ifpdf
\DeclareGraphicsExtensions{.pdf, .jpg}
\else
\DeclareGraphicsExtensions{.eps, .jpg}
\fi

\maketitle

\newpage

\begin{abstract}
We consider the prediction problem of a continuous-time stochastic process on
an entire time-interval in terms of its recent past. The approach we adopt is
based on functional kernel nonparametric regression estimation techniques where
observations are segments of the observed process considered as curves. These
curves are assumed to lie within a space of possibly inhomogeneous functions,
and the discretized times series dataset consists of a relatively small,
compared to the number of segments, number of measurements made at regular
times. We thus consider only the case where an asymptotically non-increasing
number of measurements is available for each portion of the times series. We
estimate conditional expectations using appropriate wavelet decompositions of
the segmented sample paths. A notion of similarity, based on wavelet
decompositions, is used in order to calibrate the prediction. Asymptotic
properties when the number of segments grows to infinity are investigated under
mild conditions, and a nonparametric resampling procedure is used to generate,
in a flexible way, valid asymptotic pointwise confidence intervals for the
predicted trajectories. We illustrate the usefulness of the proposed functional
wavelet-kernel methodology in finite sample situations by means of three
real-life datasets that were collected from different arenas.

\bigskip

{\em Some key words:} {\sc $\alpha$-mixing; Besov Spaces; Continuous-Time
Prediction; Functional Kernel Regression; Pointwise Prediction Intervals;
Resampling;  SARIMA Models; Smoothing Splines; Wavelets}

\end{abstract}

\section{{\sc INTRODUCTION}}
\label{sec:intro}

In many real life situations one seeks information on the
evolution of a (real-valued) continuous-time stochastic process $X
= (X(t); \; t \in \reals)$ in the future. Given a trajectory of
$X$ observed on the interval $[0, T]$, one would like to predict
the behavior of $X$ on the entire interval $[T,T+\delta]$, where
$\delta > 0$, rather than at specific time-points. An appropriate
approach to this problem is to divide the interval $[0, T]$ into
subintervals $[l\delta,(l+1)\delta]$, $l=0,1,\ldots,i-1$ with
$\delta=T/i$, and to consider the stochastic process $Z = (Z_i; \;
i \in \NN^{+})$, where $\NN^{+}=\{1,2,\ldots\}$, defined by
\begin{equation}
Z_i(t) = X(t+ (i-1)\delta), \quad i \in \NN^{+}, \quad \forall \;
\; t \in [0, \delta). \label{eqn:zprocess}
\end{equation}
This representation is especially fruitful if $X$ has a seasonal component with
period $\delta$ and can be decomposed into locally stationary parts. It can be
also employed if the data are collected as curves indexed by time-intervals of
equal lengths; these intervals may be {\em adjacent}, {\em disjoint} or even
{\em overlapping} (see, for example, Ramsay \& Silverman, 1997).

In the recent literature, practically all investigations to date
for this prediction problem are for the case, where one assumes
that an appropriately centered version of the stochastic process
$Z$ is a (zero-mean) Hilbert-valued {\em autoregressive (of order
1) processes} (ARH(1)); the best prediction of $Z_{n+1}$ given its
past history $(Z_n,Z_{n-1},\ldots,Z_1)$ is then given by
\begin{eqnarray}
\tilde{Z}_{n+1} & = & \EE(Z_{n+1} \mid Z_n, Z_{n-1},\ldots,Z_1) \nonumber\\
& = & \rho(Z_n), \quad n \in \NN^{+}, \nonumber
\end{eqnarray}
where $\rho$ is a bounded linear operator associated with the
ARH(1) process. The adopted approaches mainly differ in the way of
estimating the `prediction' operator $\rho$, or its value
$\rho(Z_n)$ given $Z_1,Z_2,\ldots,Z_n$ (see, e.g., Bosq, 1991;
Besse \& Cardot, 1996; Bosq, 2000; Antoniadis \& Sapatinas, 2003).

In many practical situations, however, the stochastic process $Z$ may not have
smooth sample paths (lying in $H$) or may not be modelled with such an
autoregressive structure. This is the case that we consider in the following
development. In particular, we also assume that the (real-valued)
continuous-time stochastic process $X = (X(t); \; t \in \reals)$ possesses a
representation of the form (\ref{eqn:zprocess}) with short duration `blocks'
$Z_i$, for $i \in \NN^{+}$. We then develop a version of prediction via
functional regression analysis, in which both the predictor and response
variables are functions of time, using a conditioning idea. Under mild
assumptions on the observed time series, a one time-interval ahead prediction
of the `block' $Z_{n+1}$ is obtained by kernel regression of the present
`block' $Z_n$ on the past `blocks' $\{Z_{n-1}, Z_{n-2}, \ldots, Z_1\}$. The
resulting predictor will be seen as a weighted average of the past `blocks',
placing more weight on `blocks' that are similar to the present one. Hence, the
analysis is rooted in the ability to find `similar blocks'. Considering that
`blocks' can be quite irregular curves, similarity matching is based on a
distance metric on the wavelet coefficients of an appropriate wavelet
decomposition of the `blocks'. A resampling scheme, involving resampling of the
original `blocks' to form `pseudo-blocks' of the same duration, is then  used
to calculate pointwise prediction intervals for the predicted `block'.

The paper is organized as follows. In Section 2, we introduce basic notions for
continuous-time prediction, relevant notions on wavelet-based orthogonal
expansions of continuous-time stochastic processes, and describe the strictly
stationarity and $\alpha$-mixing assumptions that are going to be adopted for
forecasting. In Section 3, we discuss the extension of the conditioning
approach to the one time-interval ahead prediction. Resampling-based pointwise
prediction intervals are presented in Section 4. In Section 5, we illustrate
the usefulness of the proposed functional wavelet-kernel approach for
continuous-time prediction by means of three real-life datasets that were
collected from different arenas. We also compare the resulting predictions with
those obtained by two other methods in the literature, in particular with a
smoothing spline method and with the SARIMA model. Some concluding remarks are
made in Section~6. Proofs and auxiliary results are compiled in the Appendix.

\section{{\sc PRELIMINARIES AND NOTATIONS}}
\label{sec:prelim}

\subsection{Generalities}

Let $(\Omega,\mathcal{A},\PP)$ be a probability space, rich enough so that all
random variables considered in the following development can be defined on this
space, and let $X = (X(t); \; t \in \reals)$ be a (real-valued) continuous-time
stochastic process defined on this space. Motivated by applications to
prediction and forecasting, it is supposed that the time-interval on which the
continuous-time stochastic process is observed is divided into intervals of
constant-width $\delta>0$ so that, from $X$, a functional-valued random
variable sequence $(Z_{i};~ i \in \NN^{+})$ is constructed according to the
representation (\ref{eqn:zprocess}), i.e.,
$$
Z_{i}(t)=X(t+(i-1)\delta),  \quad i \in \NN^{+}, \quad \forall \;
t\in [0,\delta).
$$
In the sequel, we regard the $Z_{i}$'s as elements of a certain (semi-)normed
functional linear space $H$ equipped with (semi-)norm $\|\cdot\|$ and its Borel
$\sigma$-field $\mathcal{F}_{\| \cdot\|}$. Recall that our aim is a one-ahead
time interval prediction which, under the above notation, is reduced in
studying a corresponding problem for the $H$-valued time series $Z=(Z_{i};~i\in
\NN^{+})$. In what follows, we assume that the time series $Z$ is strictly
stationary with $\EE(\|Z_{i}\|)<\infty$. If the time series $Z$ is not
stationary, it is assumed that it has been transformed to a stationary one by a
preprocessing procedure, and the procedures that we are going to develop hold
for the resulting stationary time series.

In practice, the random curves $Z_{i}$ are only known at discretized
equidistant time values, say $t_1, \ldots, t_P$. Thus, they must be
approximated by some $H$-valued functions, which in our case will be realized
by first expanding the $Z_{i}$'s into a wavelet basis and then estimating
consistently the coefficients from the observed discrete data. We may have used
a fixed spline basis or a Fourier basis instead but there are some good reasons
to prefer wavelet bases. When $P$ is fixed, using spline interpolation could
make sense if the sample paths exhibits a uniformly smooth temporal structure
thus not requiring any smoothing to stabilize the variance. When $P$ is large,
a necessary smoothing step in the spline approximation would be necessary and
since the smoothing is global this would not be appropriate when the sample
paths are composed of different temporal structures. The same remarks apply for
a Fourier basis. On the contrary, when $P$ is fixed and the sample paths are
continuous, one may choose an interpolation wavelet basis and the wavelets
coefficients are computed directly from the sampled values instead of inner
product integrals. When $Z_{i}$ are observed either continuously or on a very
fine discretization grid, then wavelets can be used successfully for
compression of a continuous-time stochastic process, in the sense that the
sample paths can be accurately reconstructed from a fraction of the full set of
wavelet coefficients. Whatever the setting is, the wavelet decomposition of the
sample paths will be a local one, so that if the information relevant to our
prediction problem is contained in a particular part or parts of the sample
path, as it is typically the case in many practical applications, this
information will be carried by a very small number of wavelet coefficients.

Below, we recall some background on interpolating wavelets and orthonormal
wavelet expansions of continuous-time stochastic processes that we are going to
use in the subsequent development.

\subsection{Orthonormal wavelet expansions of continuous-time stochastic processes}
\label{subsec:orthsp}

The discrete wavelet transform (DWT), as formulated by Mallat (1989) and
Daubechies (1992), is an increasingly popular tool for the statistical analysis
of time series (see, e.g., Nason \& von Sachs, 1999; Percival \& Walden, 2000;
Fryzlewicz, Van Bellegem \& von Sachs, 2003). The DWT maps a time series into a
set of wavelet coefficients, each one associated with a particular scale. Two
distinct wavelet coefficients can be either `within-scale' (i.e., both are
associated with the same scale) or `between-scale' (i.e., each has a distinct
scale).

One reason for the popularity of the DWT in times series analysis is that
measured data from most processes are inherently multiscale in nature, owing to
contributions from events occurring at different locations and with different
localization in time and frequency. Consequently, data analysis and modelling
methods that represent the measured variables at multiple scales are better
suited for extracting information from measured data than methods that
represent the variables at a single scale.

Let $S(t) \equiv S(\omega, t)$ be a mean-square continuous-time
stochastic process defined on $\Omega \times [0,1]$, i.e.
$$
S \in L^2(\Omega \times [0,1]) = \left\{ X(t): \Omega \rightarrow
\RR,\; t \in [0,1] \; \left| \; \EE \left( \int_0^1 X^2(t) \; dt
\right) < \infty \right\}. \right.
$$
Recall that (see, e.g., Neveu, 1968) $L^2(\Omega \times [0,1])$ is a separable
Hilbert space equipped with inner product defined by
$$
\langle X_1, X_2 \rangle = \EE \left( \int_0^1 X_1(t) X_2(t)\; dt
\right).
$$
To develop a wavelet-based orthonormal expansion, we mimic the construction of
a (regular) multiresolution analysis of $L^2([0,1])$ (see, e.g., Mallat, 1989).
In other words, consider a ladder of closed subspaces
$$
V_{j_0} \subset V_{j_0 +1 } \subset \cdots \subset L^2([0,1]),
$$
with any fixed $j_0\geq 0$, whose union is dense in $L^2([0,1])$, and where,
for each $j$, $V_j$ is spanned by $2^j$ orthonormal scaling functions
$\{\phi_{j,k}:~ k=0,\ldots, 2^j-1\}$, such that $\hbox{supp}(\phi_{j,k})
\subset [2^{-j}(k-c),2^{-j}(k+c)]$, with $c$ a constant not depending on $j$.
At each resolution level $j$, the orthonormal complement $W_j$ between $V_j$
and $V_{j+1}$ is generated by $2^j$ orthonormal wavelets $\{\psi_{j,k}:~
k=0,\ldots, 2^j-1\}$, such that $\hbox{supp}(\psi_{j,k}) \subset
[2^{-j}(k-c),2^{-j}(k+c)]$. As a consequence, the family $\cup_{j\geq j_0}
\{\psi_{j,k}:~k=0,\ldots,2^j-1\}$, completed with $\{
\phi_{j_0,k}:~k=0,\ldots,2^{j_0}-1\}$, constitutes an orthonormal basis of
$L^2([0,1])$.

Similarly, we define a sequence of approximating spaces of $L^2(\Omega \times
[0,1])$ by
$$
V_j(\Omega \times [0,1]) = \left\{ X \in L^2(\Omega \times [0,1]) \; \left|  \;
X(t)=\sum_{k=0}^{2^j-1} \xi_{j,k} \phi_{j,k}(t),\
\sum_{k=0}^{2^j-1}\EE(\xi_{j,k})^2 < \infty \right\}, \right.
$$
where $\{\xi_{j,k}:~k=0,\ldots,2^j-1\}$ is a sequence of random variables and
$\phi_{j,k}$ is the scaling basis of $V_j$. Note that since $L^2(\Omega \times
[0,1])$ is isomorphic to the Hilbert tensor product $L^2(\Omega) \otimes
L^2([0,1])$, the stochastic approximating spaces $V_j(\Omega \times [0,1])$ are
closed subspaces of $L^2(\Omega \times [0,1])$. Note also that for every $X \in
V_j(\Omega \times [0,1])$, one has $\EE(X) \in L^2([0,1])$ since
$$
\int_0^1\left[\EE(X(t))\right]^2 dt = \int_0^1
\left[\EE\left(\sum_{k=0}^{2^j-1} \xi_{j,k}\phi_{j,k}(t)\right)\right]^2 dt =
\sum_{k=0}^{2^j-1} \left[\EE(\xi_{j,k})\right]^2 \leq \sum_{k=0}^{2^j-1}
\EE(\xi_{j,k})^2,
$$
by orthonormality of the scaling functions. Following Cohen \& D'Ales (1997),
it is easy to see that $\{V_j(\Omega \times [0,1]):~ j \in \NN_0\}$ is a
multiresolution analysis of $L^2(\Omega \times [0,1])$. Moreover if $W_j(\Omega
\times [0,1])$ denotes the orthonormal complement of $V_j(\Omega \times [0,1])$
in $V_{j+1}(\Omega \times [0,1])$, then one naturally has the following
stochastic wavelet orthonormal expansion
$$
X \in L^2(\Omega \times [0,1]) \leftrightarrow X(t) = \sum_{k=0}^{2^{j_0}-1}
\xi_{j_0,k} \phi_{j_0,k}(t) + \sum_{j=j_0}^{\infty} \sum_{k=0}^{2^j-1}
\eta_{j,k} \psi_{j,k}(t),
$$
where $\xi_{j,k} = \int_0^1 \phi_{j,k}(t) X(t)\, dt$ and $\eta_{j,k} = \int_0^1
\psi_{j,k}(t) X(t)\, dt$. The above remarks clearly show that the wavelet-based
orthonormal expansion is a fundamental tool for viewing the continuous-time
stochastic process in both time and scale domains.

In order to allow for inhomogeneous sample paths, the notion of the Besov space
($B_{p,q}^{s}$) comes naturally into the picture. Without getting into
mathematical details, we just point out that Besov spaces are known to have
exceptional expressive power: for particular choices of the parameters $s$, $p$
and $q$, they include, e.g., the traditional H\"older ($p=q=\infty$) and
Sobolev ($p=q=2$) classes of smooth functions, and the inhomogeneous functions
of bounded variation sandwiched between $B^{1}_{1,\infty}$ and $B^{1}_{1,1}$.
The parameter $p$ can be viewed as a degree of function's inhomogeneity while
$s$ is a measure of its smoothness. Roughly speaking, the (not necessarily
integer) parameter $s$ indicates the number of function's (fractional)
derivatives, where their existence is required in an $L^p$-sense, while the
additional parameter $q$ provides a further finer gradation (see, e.g.,. Meyer,
1992). Therefore, from an approximation perspective, if the sample paths of the
continuous-time stochastic process $X$ belong to an inhomogeneous Besov space
of regularity $s > 0$, and if one uses regular enough scaling functions, one
may approximate in $L^{2}(\Omega \times [0,1])$ any sample path of the process
$X$ by its projection onto $V_J$ at a rate of the order $\mathcal{O}(2^{-s J})$
(see Cohen \& D'Ales, 1997, Theorem 2.1). This is, in fact, a simple
rephrasing, in the stochastic framework, of the deterministic results on the
multiresolution approximation of functions in Besov spaces when the error is
measured in the $L^2([0,1])$-norm. This result has the advantage that dimension
reduction by basis truncation in the wavelet domain is controlled more
precisely.

Consider now the case where the stochastic signal is sampled over a finite
number $P$ of equidistant points and assume that $P=2^J$. One then may use an
interpolating wavelet basis, as the one constructed by Donoho (1992), to
interpolate the observed values. The interpolating scaling function $\phi$
corresponds to the autocorrelation function of an orthogonal, regular enough,
scaling function and the projection onto the scaling approximating space $V_J$
is then given by
$$
{\cal P}_J (Z) (t) = \sum_{k=0}^{2^J-1} Z(t_k) \phi(2^Jt-k).
$$
While this projector has no orthogonality property, it retains however the good
approximations properties of projection operators derived from orthogonal
multiresolution analyses (see Mallat, 1999, Theorem 7.21). When $P$ is not
anymore a power of two, one may still use the above scheme by adapting the
interpolating wavelet basis to the sampling grid using an appropriate
subdivision scheme for interpolation (see Cohen, Dyn \& Matei, 2003).

We conclude this section by recalling the strictly stationarity and
$\alpha$-mixing concepts that we are going to adopt for developing the proposed
functional wavelet-kernel prediction methodology.

\subsection{Strictly stationarity and $\alpha$-mixing}
\label{subsec:orthsp}

One of our main assumption in predicting the times series $Z$ will be its
strictly stationarity. We therefore recall some results on strictly
stationarity from the above stochastic multiresolution analysis perspective.
Recall that, for all $X\in V_{J}(\Omega \times[0,1])$, we have
$$X(t)=\sum_{k=0}^{2^J-1} \xi_{J,k} \phi_{J,k}(t).$$
Therefore,
\begin{eqnarray*}
X(t+s) & =  & \sum_{k=0}^{2^J-1} \xi_{J,k} \phi_{J,k}(t+s) =
\sum_{k=0}^{2^J-1} \xi_{J,k} 2^{J/2}\phi(2^{J}(t+s) -k)\\
& = & \sum_{k=0}^{2^J-1} \xi_{J,k} 2^{J/2}\phi(2^{J}t-(k -2^{J}s)) =
 \sum_{k=0}^{2^J-1} \xi_{J,k+2^{J}s} \phi_{J,k}(t).
 \end{eqnarray*}
>From the above, it is easy to see that if $X$ is a strictly stationary process
then, at any resolution level $j$, the vector of its scaling coefficients, is
also strictly stationary. As shown in Cheng \& Tong (1996), the converse is
also true. It follows that strictly stationarity of the discrete-time series
$Z$ implies strictly stationarity in the above sense of the sequence of the
scaling coefficients vectors at any resolution. Note, moreover, that if the
strictly stationarity assumption is too strong, one could calibrate the
non-stationarity by considering only $J$-stationarity (see Cheng \& Tong,
1998), that is, strictly stationarity of the scaling coefficients up to
(finest) resolution $J$, with eventually a different distribution at each
resolution level $j$.

Our theoretical results will be derived under $\alpha$-mixing
assumptions on the time series $Z=(Z_{i};~ i \in \NN^{+})$. Recall
that for a strictly stationary series $Z=(Z_{i};~ i \in \NN^{+})$,
the $\alpha$-mixing coefficient (see Rosenblatt, 1956) is defined
by
$$\alpha_{Z}(m) = \sup_{A\in {\cal D}_{l},B\in {\cal D}_{l+m}}| \PP(A\cap B) - \PP(A) \PP(B) |,$$
where ${\cal D}_{l}=\sigma(Z_{i}, i \leq l)$ and ${\cal D}_{l+m}=\sigma(Z_{i},
i \geq l+m )$ are the $\sigma$-fields generated by $(Z_{i};~ i\leq l)$ and
$(Z_{i};~ i\geq l+m)$ respectively, for any $m \geq 1$. The stationary sequence
$Z=(Z_{i};~ i \in \NN^{+})$ is said to be $\alpha$-mixing if $\alpha_{Z}(m)
\rightarrow 0$ as $m\rightarrow \infty$.  Among various mixing conditions used
in the literature, $\alpha$-mixing is reasonably weak (see, e.g., Doukhan,
1994).

Since in the subsequent development we are dealing with wavelet decompositions,
for each $i \in \NN^{+}$, denote by $\Xi_{i}=
\{\xi_{i}^{(J,k)}:~k=0,1,\ldots,2^{J}-1\}$ the set of scaling coefficients up
to resolution level $J$ of the $i$-th segment $Z_i$. Note that because
$Z=(Z_{i};~ i \in \NN^{+})$ is a strictly stationary stochastic process, the
same is also true for the $2^{J}$-dimensional stochastic process $ \{\Xi_{i};~
i \in \NN^{+}\}$. Furthermore, denote by ${\cal A}_{J,l}=\sigma(\xi_{i}^{J,k},
i \leq l)$ and ${\cal A}_{J,l+m}=\sigma(\xi_{i}^{J,k}, i \geq l+m )$ the
$\sigma$-fields generated by $(\xi_{i}^{(J,k)};~ i\leq l)$ and
$(\xi_{i}^{(J,k)};~ i\geq l+m)$ respectively. Because $ \sigma(\xi_{i}^{(J,k)},
i \in I) \subseteq \sigma(Z_{i}, i \in I)$ for any $ I \subset \NN^{+}$, we get
\begin{eqnarray*}
\alpha_{J,k}(m) & = & \sup_{A\in{\cal A}_{J,l} ,B\in {\cal
A}_{J,l+m}} | \PP(A\cap B) - \PP(A) \PP(B) | \\ & \leq &
\sup_{A\in{\cal D}_{J,l} ,B\in {\cal D}_{J,l+m}} | \PP(A\cap B) -
\PP(A) \PP(B) | \\ & = & \alpha_{Z}(m).
\end{eqnarray*}
Note that the above observation is also true when dealing with sample paths
discretized over a fixed equidistant grid on $[0,\delta]$ of size $P$, since
then $\xi_{i}^{(J,k)}=Z_i(t_k)$ for all $k=1,2,\ldots,P$.

\section{\sc A Functional Wavelet-kernel prediction}
\label{sec:fwkp}

\subsection{Finite-dimensional Kernel Prediction}

Consider the nonparametric prediction of a (real-valued)
stationary discrete-time stochastic process. Let
$X_{n,(r)}=(X_{n}, X_{n-1},\dots, X_{n-r+1}) \in \RR^{r}$ be the
vector of lagged variables, and let $s$ be the forecast horizon.
It is well-known that the autoregression function plays a
predominant forecasting role in the above time series context.
Recall that the autoregression function $f$ is defined by
$$ f(\mathbf{x}) = \EE(X_{n+s} \mid X_{n,(r)} = \mathbf{x}), \quad \mathbf{x} \in \RR^{r}.$$
It is clear that the first task is to estimate $f$. The classical approach to
this problem is to find some parametric estimate of $f$. More specifically, it
is assumed that $f$ belongs to a class of functions, only depending on a finite
number of parameters to be estimated. This is the case of the very well-known
ARIMA models, widely studied in the literature (see, e.g., Box \& Jenkins,
1976; Brockwell \& Davis, 1991).

The above prediction problem can also be undertaken with a nonparametric view,
without any assumption on the functional form of $f$. This is a much more
flexible approach that only imposes regularity conditions on $f$. Nonparametric
methods for forecasting in time series can be viewed, up to a certain extent,
as a particular case of nonparametric regression estimation under dependence
(see, e.g., Bosq, 1991; Hart, 1991; H\"ardle \& Vieu, 1992; Hart, 1996). A
popular nonparametric method for such task is to use the kernel smoothing ideas
because they have good properties in real-valued regression problems both from
a theoretical and a computational point of view. The kernel estimator
$\hat{f}_{n}$ (based on $X_{1}, \dots, X_{n})$  of $f$ is defined by
$$ \hat{f}_{n}(\mathbf{x}) =  \frac{ \sum_{t=r}^{n-s}
\KK( (\mathbf{x} - X_{n,(r)})/h_{n})  X_{t+s} }{ \sum_{t=r}^{n-s}
\KK((\mathbf{x} - X_{n,(r)})/h_{n}) },$$ or $0$ if the denominator is null. In
our development, for simplicity, we consider a product kernel, i.e., for each
$\mathbf{x} = (x_1,\ldots,x_r)$, $\KK(\mathbf{x}) = \prod_{i=1}^{r} K(x_i)$;
also $h_{n}$ is a sequence of positive numbers (the bandwidths). The $s$-ahead
prediction is then simply given by $X_{n+s\,|\,n} = \hat{f}_{n}(X_{n,(r)})$.
Theoretical results show that the detailed choice of the kernel function does
not influence strongly the asymptotic behavior of prediction but the choice of
the bandwidth values are crucial for prediction accuracy (see, e.g., Bosq,
1998).

As it is readily seen, the prediction is expressed as a locally weighted mean
of past values, where the weights measure the similarity between
$(X_{t,(r)};~t=r,\dots,n-s)$ and $X_{n,(r)}$, taking into account the process
history. Let now $\|\cdot\|$ be a generic notation for a Euclidean norm. If the
kernel values decrease to zero as $\| \mathbf{x}\|$ increases, the smoothing
weights have high values when the $(X_{t,(r)})$ is close to $X_{n,(r)}$, and is
close to zero otherwise. In other words, the prediction $X_{n+s\,|\,n}$ is
obtained as a locally weighted average of future blocks of horizon $s$ in all
blocks of length $r$ in the past, weighted by similarity coefficients $w_{n,t}$
of these blocks with the current block, where
$$w_{n,t}(\mathbf{x})=  \frac{ \KK( (\mathbf{x} - X_{t,(r)})/h_{n}) }{
\sum_{m=r}^{n-s} \KK((\mathbf{x} - X_{m,(r)})/h_{n}) }.$$

\subsection{Functional Wavelet-Kernel Prediction}

Recall that, in our setting, the strictly stationary time series
$Z=(Z_{i};~ i \in \NN^{+})$ is functional-valued rather than
$\RR$-valued, i.e., each $Z_{i}$ is a random curve. In this
functional setup, and to simplify notation, we address, without
loss of generality, the prediction problem for a horizon $s=1$. We
could mimic the above kernel regression ideas, and use the
following estimate
\begin{equation}
Z_{n+1\,|\,n}(\cdot) = \sum_{m=1}^{n-1}w_{n,m}Z_{m+1}(\cdot),
\label{eq:Zseg}
\end{equation}
where the triangular-array of local weights $\{w_{n,m}:~m=1,2,\ldots,n-1; n \in
\NN^{+}\}$ increases with the closeness or similarity of the last observed path
$Z_{n}$ and the paths $Z_{m}$ in the past, in a (semi-)norm sense; this is made
more precise in (\ref{eq:XXseg}) below. The literature on this
infinite-dimension kernel regression related topic is relatively limited, to
our knowledge. Bosq \& Delecroix (1985) dealt with general kernel predictors
for Hilbert-valued stationary Markovian processes. A similar idea was applied
by Besse, Cardot \& Stephenson (1999) for ARH(1) processes in the special case
of a Sobolev space. Extending and justifying these kernel regression techniques
to infinite-dimensional stochastic processes with no specific structures (e.g.,
ARH(1) or more general Markovian processes), will require using
measure-theoretic assumptions on infinite-dimensional spaces (e.g., a
probability density function with respect to an invariant measure) thus
restricting the analysis and applicability to a small class of stochastic
processes (e.g., diffusion processes). Such kind of assumptions are more
natural in finite-dimensional spaces such as those obtained through orthonormal
wavelet decompositions, especially when the discretized sample paths of the
observed process are quite coarse. Taking advantage of these latter remarks,
our forecasting methodology relies upon a wavelet decomposition of the observed
curves, and uses the concepts of strictly stationarity and $\alpha$-mixing
within the stochastic multiresolution analysis framework discussed in Section
2. Moreover, note that using distributional assumptions such as those given in
the appendix on the wavelet coefficients  is much less restrictive than using
similar assumptions on the original process $Z$.

\medskip

To summarize, the proposed forecasting methodology is decomposed into two
phases
\begin{itemize}
\item[Ph1:] find within the past paths the ones that are `similar'
to the last observed path (this determines the weights);
\item[Ph2:] use the weights and the stochastic multiresoltion
analysis to forecast by a locally weighted averaging process as
the one described by (\ref{eq:Zseg}).
\end{itemize}

Since we are dealing with a wavelet decomposition, it is worth to isolate the
first phase (Ph1) by discussing possible ways to measure the similarity of two
curves, by means of their wavelet approximation, and then to proceed to the
second phase (Ph2), using again this wavelet approximation. The analysis of the
proposed kernel-based functional prediction method is based on finding similar
paths. Similarity is now defined in terms of a distance metric related to the
functional space in which the sample paths lie. When the space is a Besov
space, it is well-known that its norm is characterized by a weighted
$\ell_{p}$-norm of the wavelet coefficients of its elements (see, e.g., Meyer,
1992). It is therefore natural to address the similarity issue on the wavelet
decomposition of the observed sample paths. The wavelet transform is applied to
the observed sample paths, and due to the approximation properties of the
wavelet transform, only a few coefficients of the transformed data will be
used; a kind of contractive property of the wavelet transform.

Applying the DWT to each path, decomposes the temporal information of the time
series into discrete wavelet coefficients associated with both time and scale.
Discarding scales in the DWT that are associated with high-frequency
oscillations, provides a straightforward data reduction step and decreases the
computational burden. We want to use the distributional properties of the
wavelet coefficients of the observed series. Imagine first that we are given 2
observed series, and let $\theta^{(i)}_{j,k}$, $i=1,2$, be the discrete wavelet
coefficient of the DWT of each signal at scale $j$ ($j=j_0,\ldots,J-1$) and
location $k$ ($k=0,1,\ldots,2^j-1$). At each scale $j \geq j_0$, define a
measure of discrepancy in terms of a distance
$$
d_{j}\left(\boldsymbol{\theta}^{(1)},
\boldsymbol{\theta}^{(2)}\right) = \left\{ \sum_{k=0}^{2^{j}-1}
\left(\theta^{(1)}_{jk}- \theta^{(2)}_{j,k}\right)^{2}
\right\}^{1/2},
$$
which measures how effectively the two signals match at scale $j$.
To combine all scales, we then use
$$D\left(\boldsymbol{\theta}^{(1)}, \boldsymbol{\theta}^{(2)}\right) =
\sum_{j=j_0}^{J-1} 2^{-j}d_{j}\left(\boldsymbol{\theta}^{(1)},
\boldsymbol{\theta}^{(2)}\right).
$$
Such a measure of discrepancy is natural and is often used to test
the equality of two regression functions in the wavelet domain
(see, e.g., Spokoiny, 1996; Abramovich, Antoniadis, Sapatinas \&
Vidakovic, 2004).

As for the second phase (Ph2), recall that, for each $i \in
\NN^{+}$, $\Xi_{i}= \{\xi_{i}^{(J,k)}:~ k=0,1,\ldots,2^{J}-1\}$
denotes the set of scaling coefficients up to resolution level $J$
of the $i$-th segment $Z_i$. The kernel prediction of the scaling
coefficients at time $n+1$, $\Xi_{n+1\,|\,n}$, is given by
\begin{equation}
\Xi_{n+1\,|\,n} = \frac{ \sum_{m=1}^{n-1}K( D({\cal C}(\Xi_{n}),
{\cal C}(\Xi_{m}))/h_{n}) \Xi_{m+1} }{1/n + \sum_{m=1}^{n-1}
K(D({\cal C}(\Xi_{n}), {\cal C}(\Xi_{m}))/h_{n})},
\label{eq:XXseg}
\end{equation}
where the $1/n$ factor in the denominator allows expression
(\ref{eq:XXseg}) to be properly defined and does not affect its
convergence rate. Here, for simplicity, we use the notation $
D(x,y)/h_{n} = D(x/h_{n},y/h_{n})$, and ${\cal C}(\Xi_{k})$ is the
set of wavelet coefficients obtained by applying the ``pyramid
algorithm'' (see Mallat, 1989) on the set of (finest level)
scaling coefficients $\Xi_k$, for $k=1,2,\ldots,n$. This, leads to
the time-domain prediction at time $n+1$,
$$
Z^J_{n+1\,|\,n}(t) = \sum_{k=0}^{2^J-1} \xi_{n+1\,|\,n}^{(J,k)}
\phi_{J,k}(t), \quad \forall \; t\in [0,\delta),
$$
of $\EE(Z_{n+1}(\cdot)\,|\,Z_{n}(\cdot))$, where $\xi_{n+1\,|\,n}^{(J,k)}$ are
the components of the predicted scaling coefficients $\Xi_{n+1\,|\,n}$. The
following theorem shows its consistency property.

\vspace{0.2cm}

\begin{theorem} Suppose that the Assumptions (A1)-(A7), given in the Appendix, are true.
\begin{itemize}
\item[(i)] If $t \in \{ 0,1,\ldots,2^{J}-1 \} $ and if $h_n=\left(\frac{\log n}{n}\right)^{1/(2+2^J)}$, then,
${\rm as} \quad  n \rightarrow \infty$, we have
$$
\Big| Z^J_{n+1\,|\,n}(t) \ - \ \EE(Z_{n+1}(t)\,|\,Z_{n}) \Big| = O\left(2^{J/2}
\left(\frac{\log n}{n}\right)^{1/(2+2^J)} \right),\quad \textit{almost surely}.
$$
\item[(ii)] If the sample paths are sampled on a grid of size $2^J$ and if $h_n=\left(\frac{\log n}{n}\right)^{1/(2+2^J)}$,
then, ${\rm as} \quad  n \rightarrow \infty$, we have
$$
\Big| Z^J_{n+1\,|\,n}(t) \ - \ \EE(Z_{n+1}(t)\,|\,Z_{n}) \Big| = O\left(2^{J/2}
\left(\frac{\log n}{n}\right)^{1/(2+2^J)} + 2^{-sJ} \right), \quad
\textit{almost surely}.
$$
\end{itemize}
\label{th:cons}
\end{theorem}

\begin{remark}
{\rm  Note that in both assertions of Theorem~\ref{th:cons}, the size of the
sampling grid over each segment affects the convergence rate of our predictor.
In the first case, which is the most usual in practice, the rate becomes slower
as the dimension $P$ increases but we still have consistency as the number of
segments increases to infinity. In the second case, however, an extra term is
given by the wavelet approximation of the sample paths at resolution $J$ and
getting a larger $J$ affects considerably the rate of the estimator which is
the well-known problem of the ``curse of dimensionality''. One possible way to
deal with this problem would be to look at the regressor in an
infinite-dimensional space but, as already noted above, this is a difficult
problem, since one would need some concentration assumption about the
distribution of the functional-valued time series $Z=(Z_{i};~ i \in \NN^{+})$
without referring to any particular probability density function.}
\end{remark}

We conclude this section by pointing out that, as in any
nonparametric smoothing approach, the choice of the smoothing
parameter $h_n$ (the bandwidth) is of great importance. Once $h_n$
is specified, only time segments that lie within a similarity
distance from the segment $Z_n$ within $h_n$ will be used to
estimate the prediction at time $n+1$. Intuitively, a large value
of $h_n$ will lead to an estimator that incurs large bias, while a
small value, might reduce the bias but the variability of the
predicted curve could be large since only few segments are used in
the estimation. A good choice of $h_n$ should balance the
bias-variance trade off. In our implementation, we have used the
leave-one out cross-validation for times-series data as suggested
by Hart (1996). The principle of the cross-validation criterion is
to select the bandwidth which, for our given prediction horizon
$s=1$, minimizes the mean squared prediction errors of the
$(i+1)$-th segment using all segments in the past except the
$i$-th, i.e., the value of $h_n$ that minimizes
$$
CV(h) = \frac{1}{n-1}\sum_{i=1}^{n-1} \left\|Z_{i+1} -
Z^{(i)}_{i+1\,|\,i}\right\|^2,
$$
where $Z^{(i)}_{i+1\,|\,i}$ is the kernel regression estimate with bandwidth
$h$ obtained using the series without its $i$-th segment. This is the method
for choosing the bandwidth adopted in the numerical results presented in
Section~\ref{sec:num}.

\section{{\sc Resampling-based Pointwise Prediction Intervals}}
\label{sec:boot}

Apart from the prediction $Z^J_{n+1\,|\,n}(t)$ discussed in
Section~\ref{sec:fwkp}, we also construct resampling-based pointwise prediction
intervals for $ Z_{n+1}(t)$. In particular, suppose that $ Z_{n}(t)$ is
observed at the set $ 0 \leq t_{1} < t_{2} < \ldots < t_{P}\leq \delta $ of
discrete points on the interval $ [0,\delta]$. A pointwise prediction interval
for $ Z_{n+1}(t)$ is defined to be a set of lower and upper points $
L_{n+1,\alpha}(t_{i})$ and $ U_{n+1,\alpha}(t_{i})$ respectively, such that for
every $ t_{i}, i=1,2, \ldots, P$, and every $\alpha \in (0,1)$,
\[
\PP \left(L_{n+1,\alpha}(t_{i}) \leq Z_{n+1}(t_{i}) \leq
U_{n+1,\alpha}(t_{i})\right) \geq 1 - 2\alpha.
\]
Note that since we are looking at the one step prediction of $ Z_{n+1}(t)$
given $ Z_{n}$, we are in fact interested in the conditional distribution of
$Z_{n+1}(t)$ given $ Z_{n}$, i.e., $L_{n+1,\alpha}(t_{i})$ and $
U_{n+1,\alpha}(t_{i}) $ are the lower and upper $\alpha$-percentage points of
the conditional distribution of $ Z_{n+1}(t_{i})$ given $ Z_{n}$.

To construct such a prediction interval the following simple resampling
procedure is proposed. Given $ Z_{n}$, i.e., given
 $ C(\Xi_{n})$, define the weights
\[ w_{n,m} = \frac{ K( D({\cal C}(\Xi_{n}),
{\cal C}(\Xi_{m}))/h_{n})  }{n^{-1} + \sum_{m=1}^{n-1} K(D({\cal C}(\Xi_{n}),
{\cal C}(\Xi_{m}))/h_{n})} + \frac{1}{1+ n\sum_{m=1}^{n-1} K(D({\cal
C}(\Xi_{n}), {\cal C}(\Xi_{m}))/h_{n}))}.\] Note that the weights have been
selected appropriately so that
\[ 0 \leq w_{n,m} \leq 1 \ \ \ \ \mbox{and} \ \ \ \ \sum_{m=1}^{n-1}w_{n,m}=1.\]
Now, given $ Z_{n}$, generate pseudo-realizations $ Z_{n+1}^{\ast}(t)$  such
that for $ m=1,2, \ldots, n-1$,
\[
\PP\left(Z_{n+1}^{\ast}(t) = Z_{m+1}(t) \,|\,  Z_{n} \right) = w_{n,m},
\]
i.e., $ Z_{n+1}^{\ast}(t)$ is generated by choosing randomly a segment from the
whole  set of observed segments $ Z_{m+1}(t)$, where the probability that the
$(m+1)$-th segment is chosen depends on how `similar' is the preceding segment
$ Z_{m}$ to $ Z_{n}$. This `similarity' is measured by the resampling
probability $ w_{n,m}$.

Given pseudo-replicates $ Z_{n+1}^{\ast}(t)$, calculate $ R_{n+1}^{\ast}(t_{i})
= Z_{n+1}^{\ast}(t_{i}) - Z_{n+1|\,n}^{J}(t_{i})$, where $
Z_{n+1|\,n}^{J}(t_{i})$ is our time-domain conditional mean predictor. Let $
R_{n+1,\alpha}^{\ast}(t_{i})$ and  $ R_{n+1,1-\alpha}^{\ast}(t_{i})$ be the
lower and upper $\alpha$-percentage points of $ R_{n+1}^{\ast}(t_{i})$. Note
that these percentage points can be consistently  estimated  by the
corresponding empirical percentage points over $B$ realizations $
Z_{n+1}^{\ast^{(b)}}(t_{i})$, $ b=1,2, \ldots, B$,  of  $ Z_{n+1}^{\ast}(t_i)$.
A $(1-\alpha)100\%$ pointwise prediction interval for $ Z_{n+1}(t_{i})$ is then
obtained by
\[ \{ [\ L_{n+1,\alpha}^{\ast}(t_{i}),\  U^{\ast}_{n+1,\alpha}(t_{i})\ ], \ \ i=1,2, \ldots, P\}, \]
where $ L^{\ast}_{n+1,\alpha}(t_{i})=  R_{n+1,\alpha}^{\ast}(t_{i})+
Z_{n+1|\,n}^{J}(t_i)$ and $
U_{n+1,\alpha}^{\ast}(t_{i})=R_{n+1,1-\alpha}^{\ast}(t_{i})+
Z_{n+1|\,n}^{J}(t_i) $.

The following theorem shows that the proposed method is asymptotically valid,
i.e., the so-constructed resampling-based prediction interval achieves the
desired pointwise coverage probability.

\vspace*{0.3cm}

\begin{theorem}
Suppose that the Assumptions (A1)-(A7), given in the Appendix, are true. Then,
for every $ i=1,2, \ldots, P$, and every $\alpha \in (0,1)$
 \label{th:boot}
\[\lim_{n\rightarrow\infty}\PP \left(L_{n+1,\alpha}^{\ast}(t_{i}) \leq Z_{n+1}(t_{i}) \leq
U_{n+1,\alpha}^{\ast}(t_{i})\,|\, Z_{1}, \ldots, Z_{n}\right) \geq 1 -
2\alpha.\]
\end{theorem}

\section{{\sc Applications to real-life datasets}}
\label{sec:num}

We now illustrate the usefulness of the proposed functional wavelet-kernel
({\sc W-K}) approach for continuous-time prediction in finite sample situations
by means of three real-life datasets that were collected from different arenas,
in particular with $(a)$ the prediction of the entire annual cycle of
climatological El Ni\~no-Southern Oscillation time series one-year ahead from
monthly recordings, $(b)$ the one-week ahead prediction of Paris electrical
load consumption from half-hour daily recordings, and $(c)$ the one-year ahead
prediction of the Nottingham temperature data from monthly recordings.

For the {\sc W-K} approach, the interpolating wavelet transform of Donoho
(1992) based on {\em Symmlet 6} (see Daubechies, 1992, p. 195) was used.
Preliminary simulations show that the analysis is robust with respect to the
wavelet filter, e.g., using {\em Coiflet 3} (see Daubechies, 1992, p. 258). In
the case where the number of time points ($P$) in each segment is not a power
of 2, each segment is extended by periodicity at the right to a length closest
to the nearest power of 2. The Gaussian kernel $(K)$ was adopted in our
analysis. Again, preliminary simulations show that our analysis is robust with
respect to kernels with unbounded support (e.g., Laplace). The bandwidth
($h_n$) was chosen by the leave-one out cross-validation for times-series data
as suggested by Hart (1996). For the associated 95\% resampling-based pointwise
prediction intervals, the number of resampling samples ($B$) was taken equal to
500.

We compare the resulting predictions with those obtained by some
well-established methods in the literature, in particular with a smoothing
spline ({\sc SS}) method and with the classical {\sc SARIMA} model. The {\sc
SS} method, introduced by Besse \& Cardot (1996), assumes an ARH(1) structure
for the time series $Z=(Z_{i};~ i \in \NN^{+})$ and handles the discretization
problem of the observed curves by simultaneously estimating the sample paths
and projecting the data on a $q$-dimensional subspace (that the predictable
part of $Z$ assumed to belong) using smoothing splines (by solving an
appropriate variational problem). The corresponding smoothing parameter
($\lambda$) and dimensionality ($q$) are chosen by a cross-validation
criterion. Following the Box-Jenkins methodology (see Box \& Jenkins, 1976,
Chapter 9), a suitable {\sc SARIMA} model is also adjusted to the times series
$Z=(Z_{i};~ i \in \NN^{+})$.

The quality of the prediction methods was measured by the {\em relative
mean-absolute error} (RMAE) defined by
\begin{equation}
\text{RMAE} = \frac{1}{P} \sum_{t=1}^{P} \frac{|\hat{Z}_{n_0}(t_i) -
Z_{n_0}(t_i)|}{|Z_{n_0}(t_i)|}, \label{eqn:rmae}
\end{equation}
where $Z_{n_0}$ is the $n_0$-th element of the time series $Z$ and
$\hat{Z}_{n_0}$ is the prediction of $Z_{n_0}$ given the past.

The computational algorithms related to wavelet analysis were
performed using Version 8.02 of the freeware WaveLab software. The
entire numerical study was carried out using the Matlab
programming environment.

\subsection{El Ni\~no-Southern Oscillation} \label{subsec:elnino}

This application concerns with the prediction of a climatological times series
describing El Ni\~no-Southern Oscillation (ENSO) during the 12-month period of
1986, from monthly observations during the 1950--1985 period. ENSO is a natural
phenomenon arising from coupled interactions between the atmosphere and the
ocean in the tropical Pacific Ocean. El Ni\~no (EN) is the ocean component of
ENSO while Southern Oscillation (SO) is the atmospheric counterpart of ENSO.
Most of the year-to-year variability in the tropics, as well as a part of the
extra-tropical variability over both Hemispheres, is related to ENSO. For a
detailed review of ENSO the reader is referred, for example, to Philander
(1990).

\figh{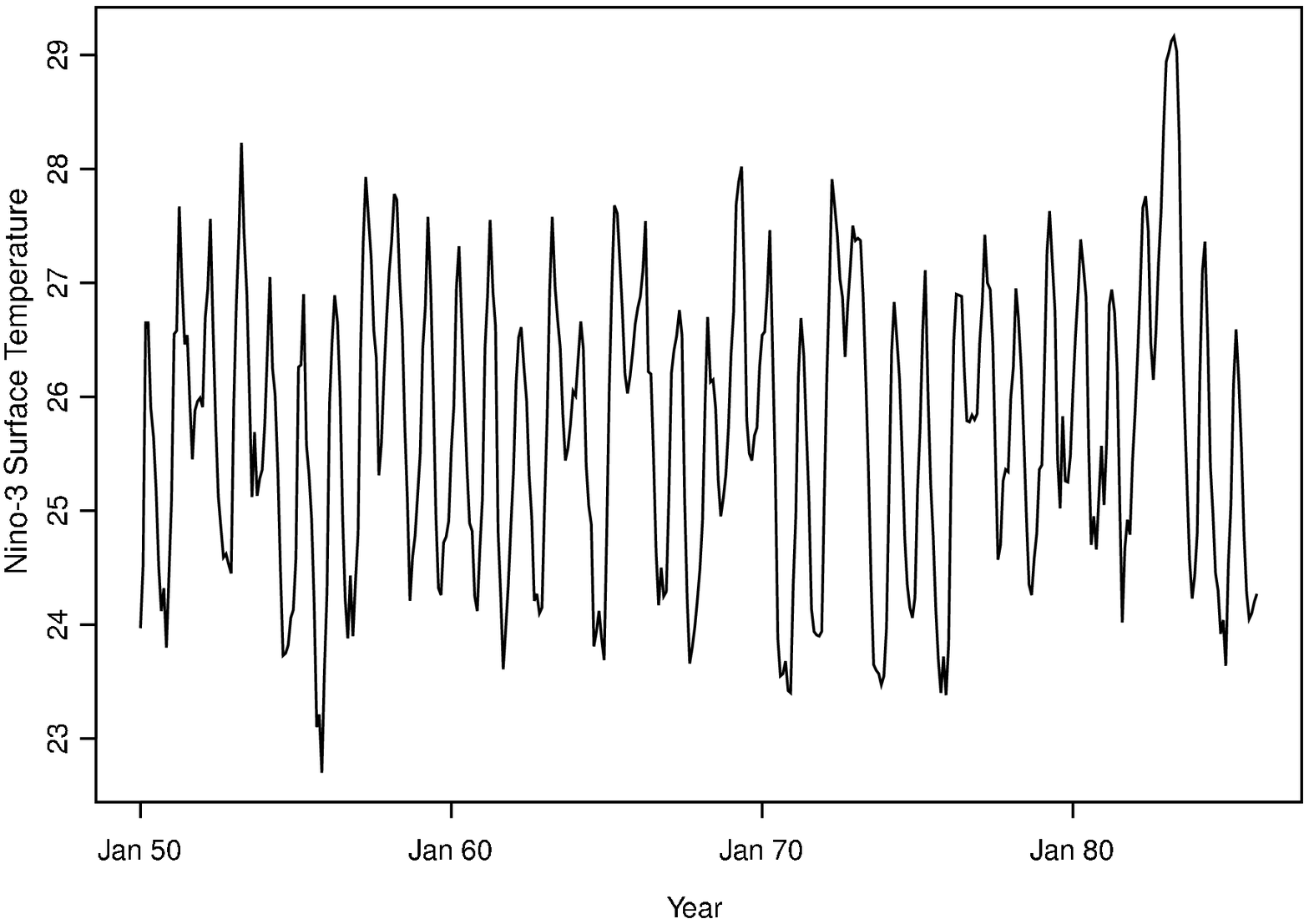}{5in}{The monthly mean Ni\~no-3 surface
temperature index in (deg C) which provides a contracted
description of ENSO.}

An useful index of El Ni\~no variability is provided by the sea surface
temperatures averaged over the Ni\~no-3 domain ($5^{\rm o}{\rm S}-5^{\rm o}{\rm
N}$, $150^{\rm o}{\rm W}-90^{\rm o}{\rm W}$). Monthly mean values have been
obtained from January 1950 to December 1996 from gridded analyses made at the
U.S. National Centers for Environmental Prediction (see Smith, Reynolds,
Livezey \& Stokes, 1996). The time series of this EN index is depicted in
Figure~\ref{elninoseries}, and shows marked inter-annual variations
superimposed on a strong seasonal component. It has been analyzed by many
authors (see, for example, Besse, Cardot \& Stephenson, 2000; Antoniadis \&
Sapatinas, 2003).

The bandwidth ($h_n$) for the {\sc W-K} method was chosen by a cross-validation
criterion and found equal to 0.11. We have compared our results with those
obtained by Besse, Cardot \& Stephenson (2000), using the {\sc SS} method, with
smoothing parameter ($\lambda$) and dimensionality ($q$) chosen optimally by a
cross-validation criterion and found equal to $1.6\times10^{-5}$ and 4,
respectively. To complete the comparison, a suitable {\sc ARIMA} model,
including 12 month seasonality, has also been adjusted to the times series from
January 1950 to December 1985, and the most parsimonious {\sc SARIMA} model,
validated through a portmanteau test for serial correlation of the fitted
residuals, was selected.

\figh{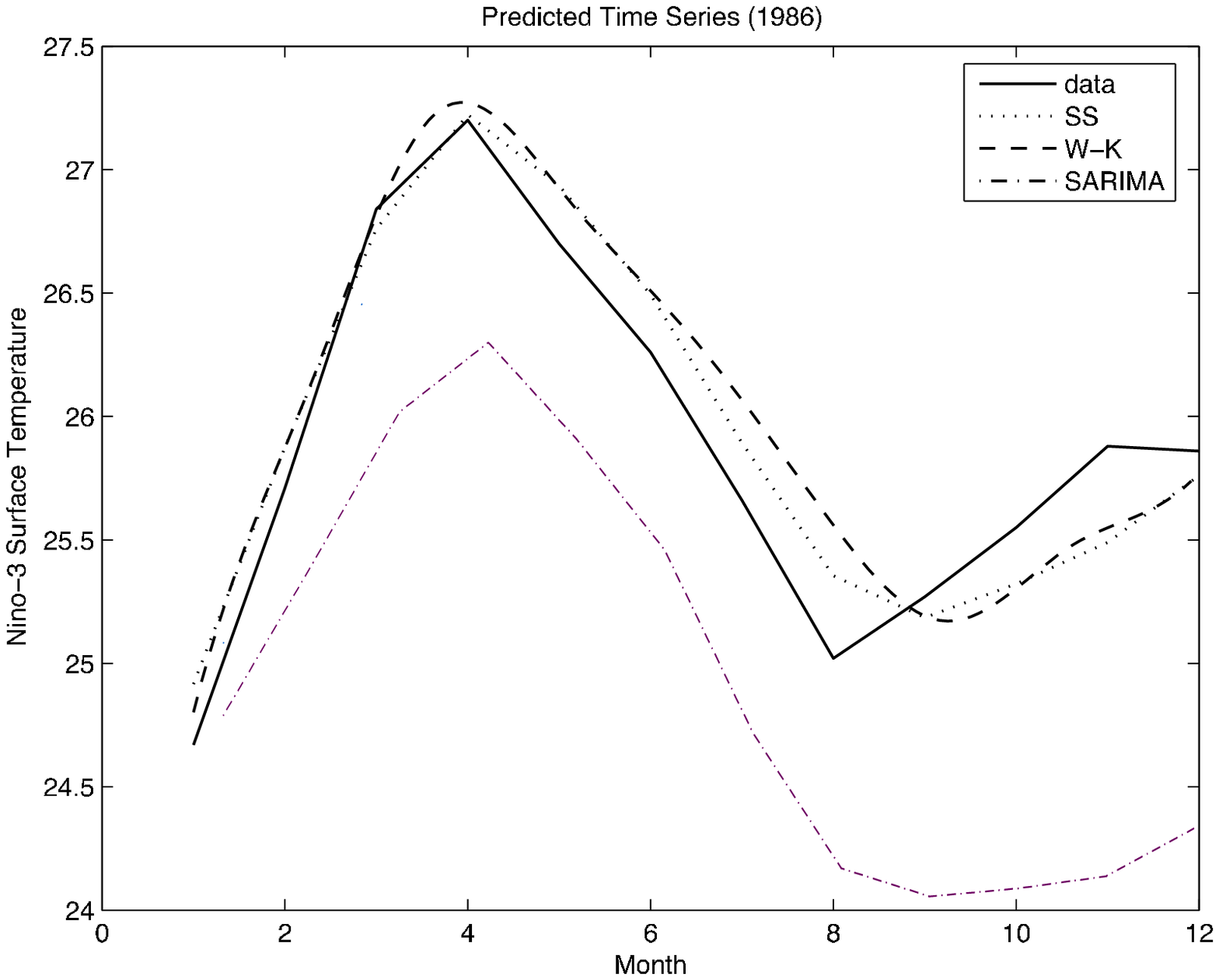}{5in}{The Ni\~no-3 surface temperature during
1986 (---) and its various predictions using the {\sc W-K}
(-\,-\,-), {\sc SS} ($\cdots$), and {\sc SARIMA} (-\,$\cdot$\,-)
methods.}

\figh{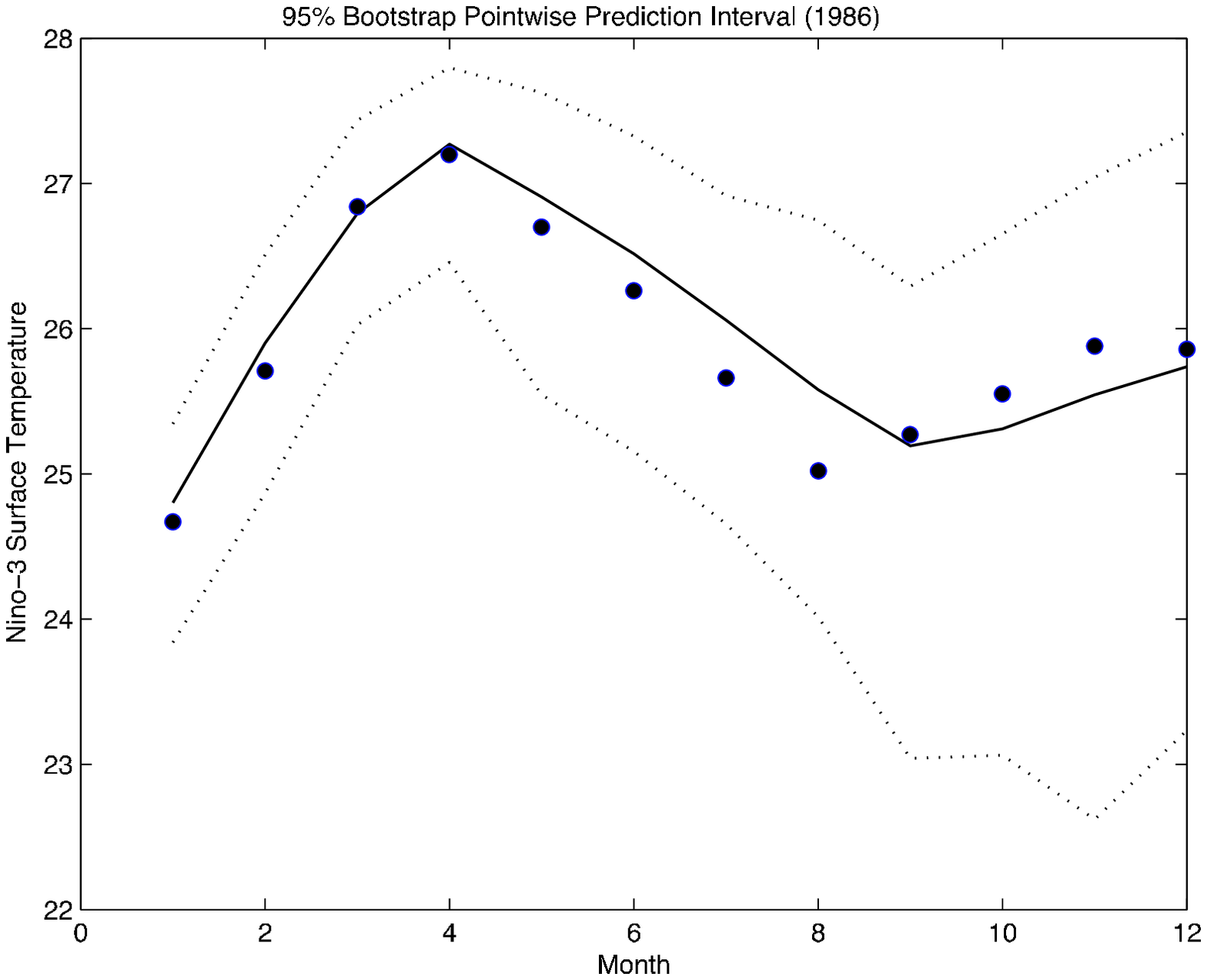}{5in}{95\% resampling-based pointwise prediction interval
($\cdots$) for the Ni\~no-3 surface temperature during 1986, based on the
corresponding prediction obtained by the {\sc W-K} (---) method. The true
points ($\bullet$) are also displayed.}

Figure~\ref{elninopred} displays the observed data of the 37th year (1986) and
its predictions obtained by the {\sc W-K}, {\sc SS} and {\sc SARIMA} methods.
The RMAE of each prediction method are displayed in Table~1 (we have taken
$n_0= 37$ and $P=12$). As observed in both the figure and the table, the {\sc
W-K} and {\sc SS} estimators give almost similar predictions, both visually and
in terms of RMAE. The prediction obtained by the {\sc SARIMA} model is strongly
and uniformly biased, failing thus to produce an adequate prediction (see
Besse, Cardot \& Stephenson, 2000, for an explanation). Note that, after May,
all predictions are not very close to the true points and this difficulty in
prediction is captured in Figure~\ref{CIelnino} which displays the
corresponding 95\% resampling-based pointwise prediction interval for the
Ni\~no-3 surface temperature during 1986, based on the corresponding prediction
obtained by the {\sc W-K} method. As observed in the figure, it becomes clear
that as one moves from May onwards, this interval gets larger.

\begin{table}
\label{errors-elnino}
\begin{center}
\begin{tabular}{|c||c|}
\hline
{\bf Prediction Method} & {\bf RMAE} \\
\hline \hline
    {\sc W-K} &  0.86 \% \\
    {\sc SS} &  0.76 \% \\
    {\sc SARIMA} & 3.72 \% \\
\hline
\end{tabular}
\end{center}
\caption{RMAE for the prediction of Ni\~no-3 surface temperatures during 1986
based on the {\sc W-K}, {\sc SS} and {\sc SARIMA} methods.}
\end{table}

\subsection{Paris Electrical Load Consumption} \label{subsec:elec}

This application concerns with the one-week ahead prediction of Paris
electrical load consumption from half-hour daily recordings. The short-term
predictions are based on data sampled over 30 minutes, obtained after
eliminating certain components linked to weather conditions, calendar effects,
outliers and known external actions. The dataset analyzed is part of a larger
series recorded from the French national electricity company (EDF) during the
period running from the 1st of August 1985 to the 4th of July 1992. The time
period that we have analyzed runs 35 days, starting from the 24th of July 1991
to the 27th of August 1991, and it is displayed in Figure~\ref{elecseries}. One
may note quite a regularity in this time series and a marked periodicity of 7
days (linked to economic rhythms) together with a pseudo-daily periodicity.
However, daily consumption patterns due to holidays, weekends and discounts in
electricity charges (e.g., relay-switched water heaters to benefit from special
night rates), make the use of {\sc SARIMA} modelling for forecasting
problematic for about 10\% of the days when working with half-hour data (see
Misiti, Misiti, Oppenheim \& Poggi, 1994).

\figh{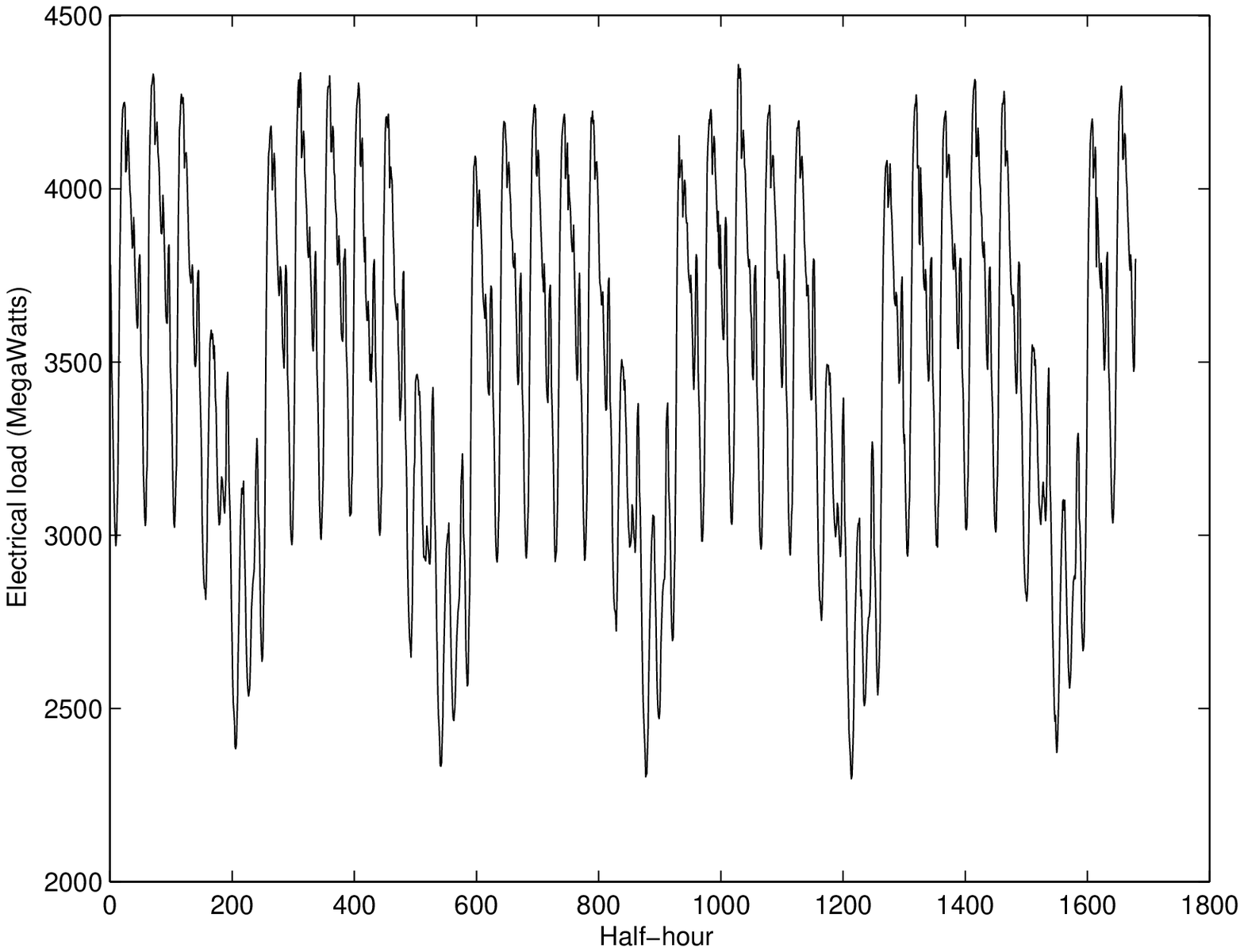}{5in}{The half-hour electricity load consumption
in Paris from the 24th of July 1991 to the 27th of August 1991.}

\figh{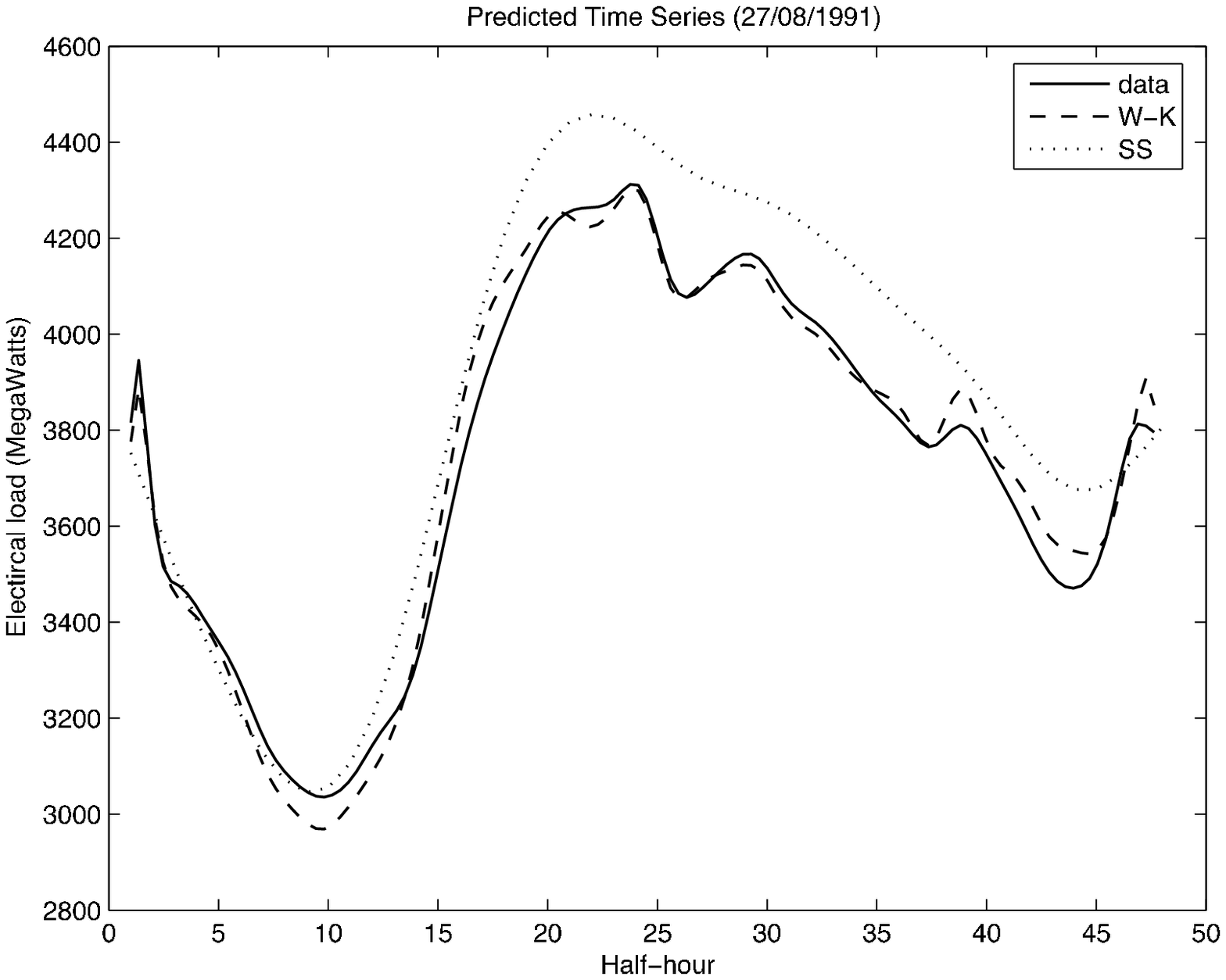}{5in}{The half-hour electricity load consumption in
Paris during 27th of August 1991 (---) and its various predictions
using the {\sc W-K} (-\,-\,-) and {\sc SS} ($\cdots$) methods.}

The bandwidth ($h_n$) for the {\sc W-K} method was chosen by cross-validation
and found equal to 0.01. We have compared our results with those obtained using
the {\sc SS} method, with smoothing parameter ($\lambda$) and dimensionality
($q$) chosen by cross-validation and found equal to $10^{-4}$ and 4,
respectively. Figure~\ref{elecpred} displays the observed data of the 2th
August 1991 and its predictions obtained only by the {\sc W-K} and {\sc SS}
methods. The RMAE of both prediction methods are displayed in Table~2 (we have
taken $n_0=35$ and $P=48$). As observed in both the figure and the table, the
prediction obtained by the {\sc W-K} method is reasonably close to the true
points, while the prediction made by the {\sc SS} method falls far off from
them. This example, clearly illustrates the impact of the proposed
wavelet-based approach since the trajectory to be predicted seems not regular
with some peculiar peaks. On the other hand, the smoothing spline-based
approach relies upon more stringent smoothness assumptions and the
corresponding prediction is therefore largely biased, falling well off the true
points. Figure~\ref{CIelec} displays the 95\% resampling-based pointwise
prediction interval for the half-hour electricity load consumption in Paris
during 27th of August 1991, based on the corresponding prediction obtained by
the {\sc W-K} method.

\figh{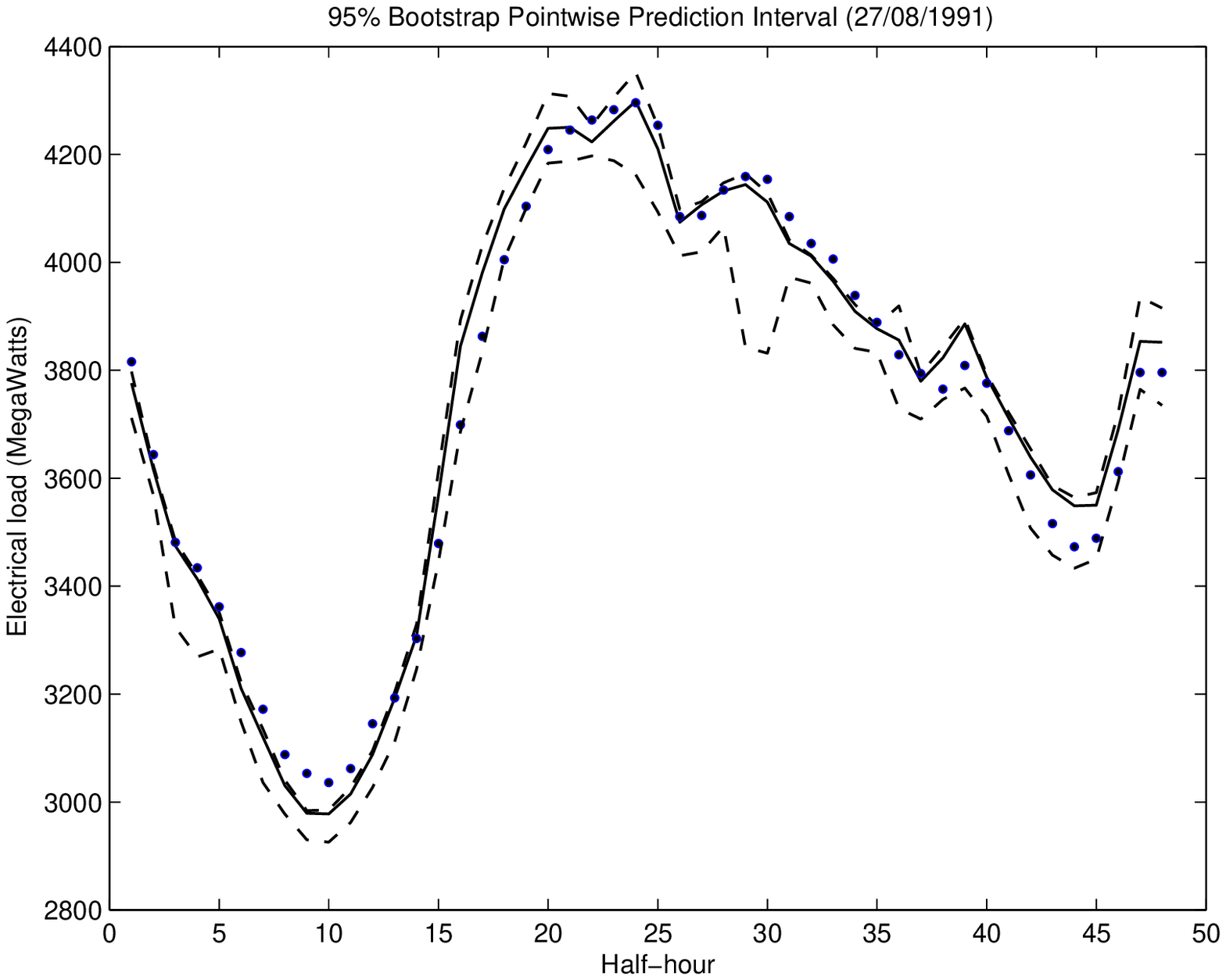}{5in}{95\% resampling-based pointwise prediction interval
(-\,-\,-) for the half-hour electricity load consumption in Paris during 27th
of August 1991, based on the corresponding prediction obtained by the {\sc W-K}
(---) method. The true points ($\bullet$) are also displayed.}

\begin{table}
\label{errors-elec}
\begin{center}

\begin{tabular}{|c||c|}
\hline
{\bf Prediction Method} & {\bf RMAE} \\
\hline \hline
    {\sc W-K} & 12 \% \\
    {\sc SS} & 36 \% \\
\hline
\end{tabular}
\end{center}

\caption{RMAE for the prediction of half-hour electricity load consumption in
Paris during 27th of August 1991 based on the {\sc W-K} and {\sc SS} methods.}
\end{table}

\subsection{Nottingham Temperature Data} \label{subsec:nottin}

This application concerns with the one-year ahead prediction of the Nottingham
temperature data from monthly recordings. The dataset analyzed are mean monthly
air temperatures ($^{\rm o}$F) at Nottingham castle from January 1920 to
December 1939, from `Meteorology of Nottingham', in {\em City Engineer
Surveyor}. Since February 1929 was an exceptionally cold month in England, and
since the seasonal pattern is fairly stable over time, the original dataset has
been `corrected'. In other words, since this `outlier' will distort the fitting
process, the value of February 1929 was altered it to a low value for February
of $35^{\rm o}$F (see Venables \& Ripley, 1999, Chapter 13). This `corrected'
dataset, which is the series {\tt nottem} in the {\tt MASS} library of {\sc
S-PLUS}, is the one analyzed below.

\figh{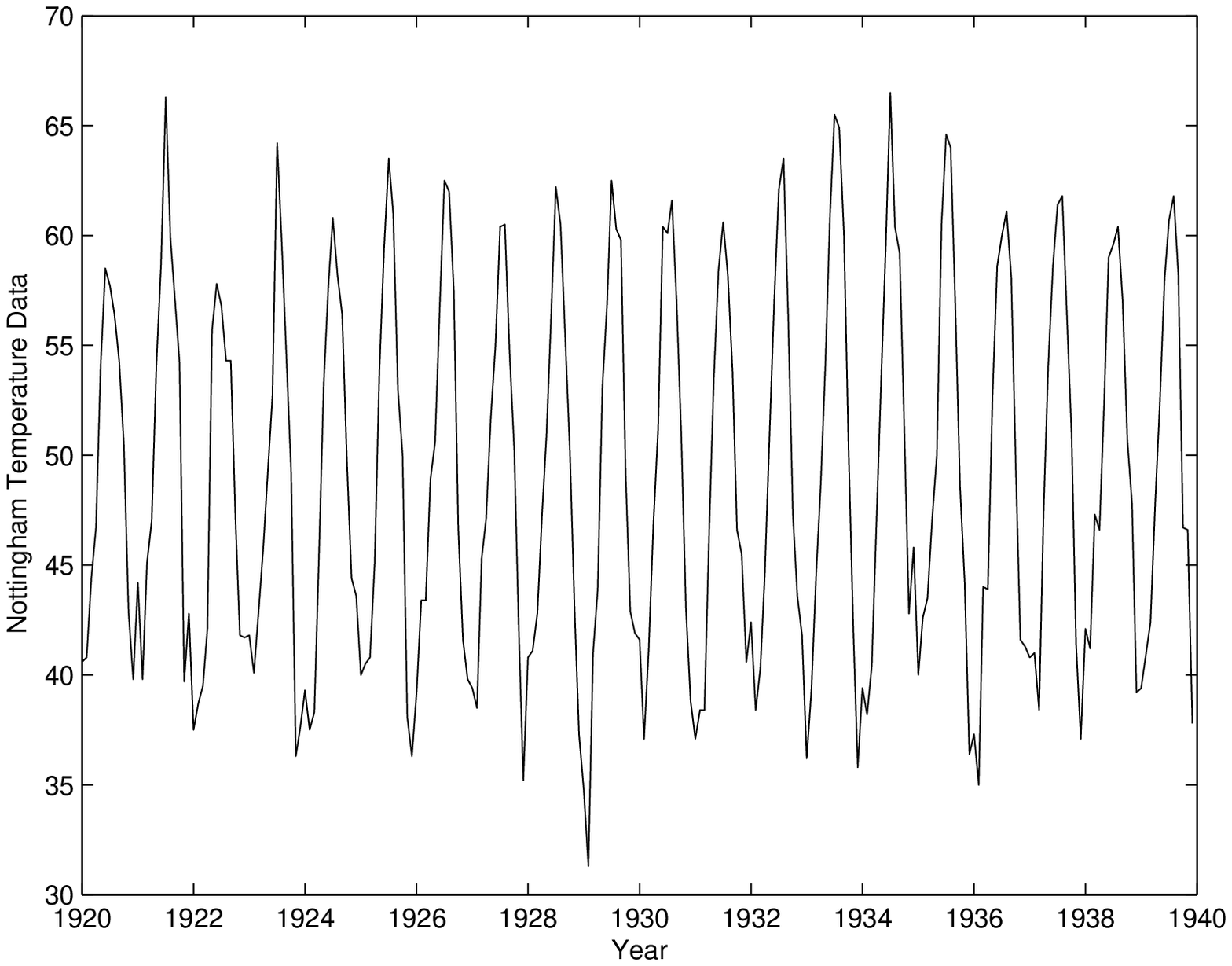}{5in}{The mean monthly air temperatures ($^{\rm
o}$F) at Nottingham castle from January 1920 to December 1939.}

\figh{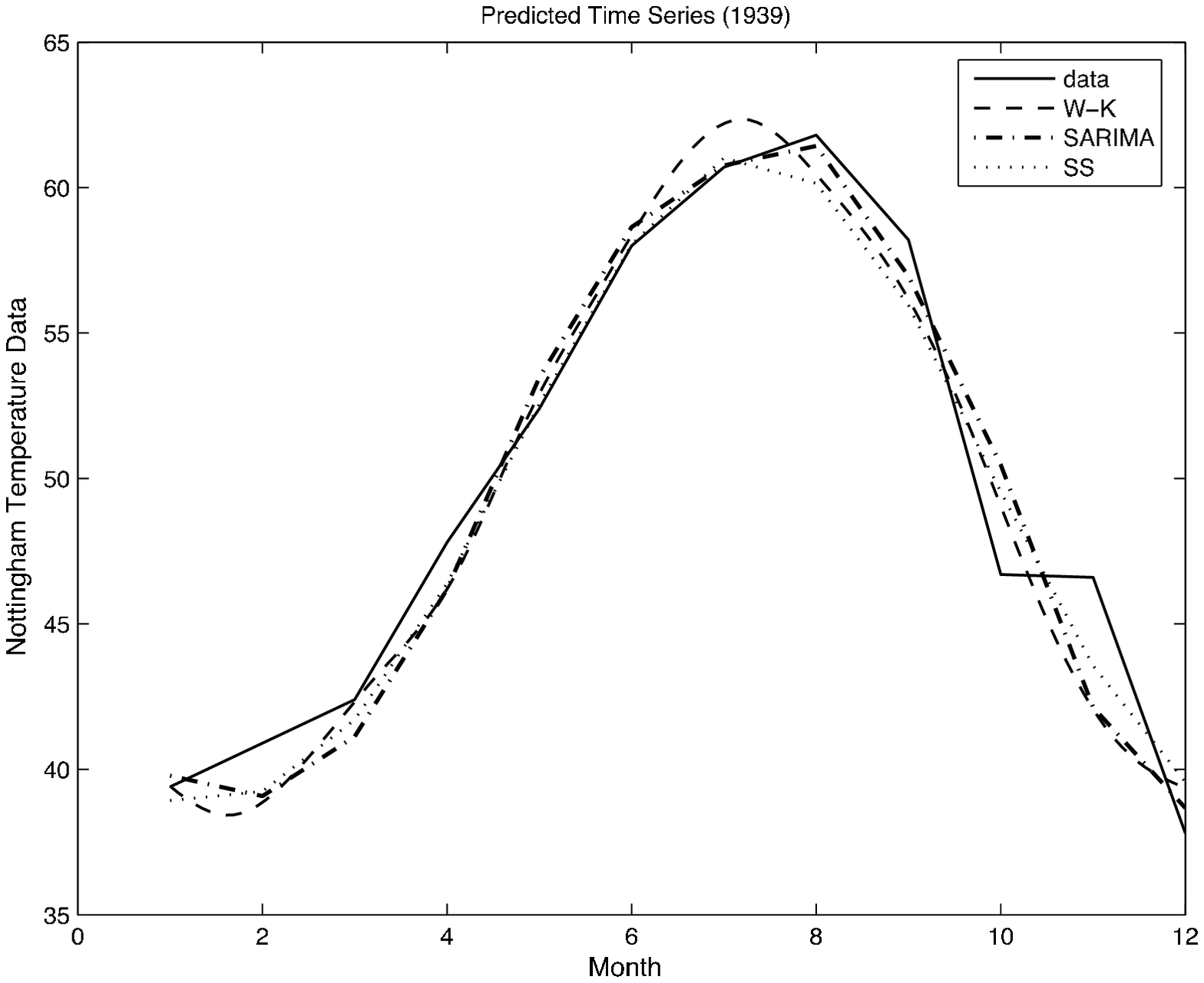}{5in}{The mean monthly air temperatures ($^{\rm
o}$F) at Nottingham castle during 1939 (---) and its various
predictions using the {\sc W-K} (-\,-\,-), {\sc SS} ($\cdots$),
and {\sc SARIMA} (-\,$\cdot$\,-) methods.}

The bandwidth ($h_n$) for the {\sc W-K} method was chosen by cross-validation
and found equal to 2.5. We have compared our results with those obtained using
the {\sc SS} method, with smoothing parameter ($\lambda$) and dimensionality
($q$) chosen by cross-validation and found equal to $10^{-4}$ and 1,
respectively. To complete the comparison, a suitable {\sc ARIMA} model,
including 12 month seasonality, has also been adjusted to the times series from
January 1920 to December 1938, and the most parsimonious {\sc SARIMA} model,
validated through a portmanteau test for serial correlation of the fitted
residuals, was selected. Figure~\ref{nottpred} displays the observed data of
the year 1939 and its various predictions obtained by the {\sc W-K}, {\sc SS}
and {\sc SARIMA} methods. It is well-known in the literature (see, e.g.,
Venables \& Ripley, 1999, Chapter 13) that this time series can be adequately
predicted by a parametric model ({\sc SARIMA}), but it is evident that both
nonparametric models ({\sc W-K} and {\sc SS}) perform equally-well. The RMAE of
each prediction method is displayed in Table~3 (we have taken $n_0=20$ and
$P=12$). As observed in both the figure and the table, before March and after
October, all predictions are largely biased and hence not very close to the
true points. This difficulty in prediction is captured in Figure~\ref{CIelec}
which displays the corresponding 95\% resampling-based pointwise prediction
interval for the mean monthly air temperatures ($^{\rm o}$F) at Nottingham
castle during 1939, based on the corresponding prediction obtained by the {\sc
W-K} method. As observed in the figure, it becomes clear that as one moves from
March backwards and from October onwards, this interval gets larger.

\figh{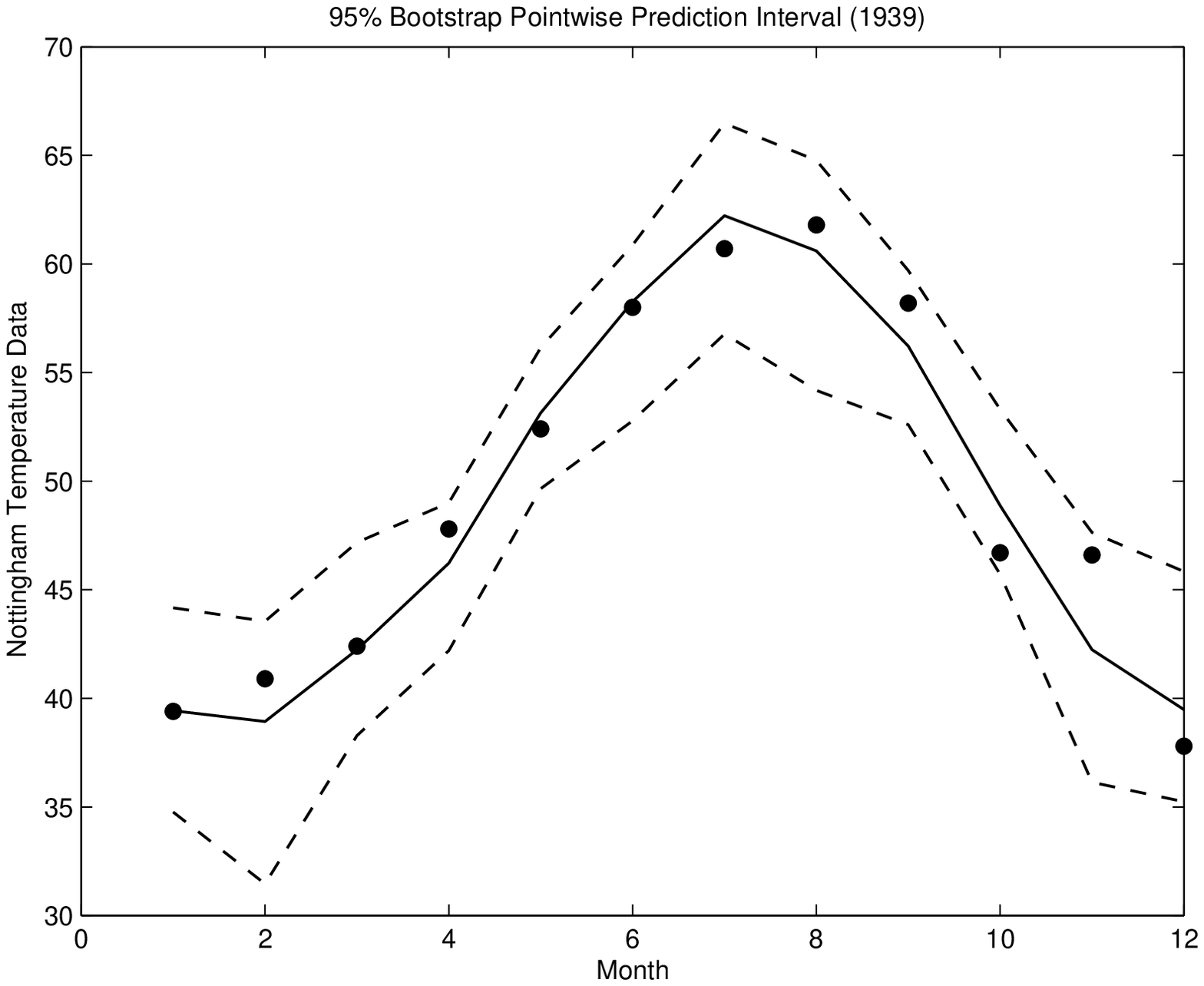}{5in}{95\% resampling-based pointwise prediction interval
(-\,-\,-) for the mean monthly air temperatures ($^{\rm o}$F) at Nottingham
castle during 1939, based on the corresponding prediction obtained by the {\sc
W-K} (---) method. The true points ($\bullet$) are also displayed.}

\begin{table}
\label{errors-pred}
\begin{center}

\begin{tabular}{|c||c|}
\hline
{\bf Prediction Method} & {\bf RMAE} \\
\hline \hline
    {\sc W-K} & 30 \% \\
    {\sc SS} & 28 \% \\
    {\sc SARIMA} & 31 \% \\
\hline
\end{tabular}
\end{center}

\caption{RMAE for the prediction of half-hour electricity load consumption in
Paris during 27th of August 1991 based on the {\sc W-K} and {\sc SS} methods.}
\end{table}

\section{{\sc Conclusions}}
The functional wavelet-kernel prediction methodology, of a continuous-time
stochastic process on an entire time-interval in terms of its recent past,
developed in this paper exhibits very good performance with respect to other
well-known parametric and nonparametric techniques. As it is demonstrated in
the real-life datasets analyzed, the proposed wavelet-based prediction
methodology outperforms the smoothing spline-based prediction methodology for
stochastic processes with inhomogeneous sample paths, and performs equally-well
for stochastic processes with quite regular sample paths. Moreover, it performs
reasonably well in situations where the classical seasonal parametric model
exhibits very good predictions. We have, however, noted than when the number of
sampling points within each time-segment is large, the curse of dimensionality
leads to inefficiency of the proposed nonparametric prediction method unless
the number of segments is really large.

\section*{Acknowledgements}
Anestis Antoniadis was supported by the `IAP Research Network
P5/24'  and the `Cyprus-France CY-FR/0204/04 Zenon Program'. Efstathios
Paparoditis and Theofanis Sapatinas were supported by the
`Cyprus-France CY-FR/0204/04 Zenon Program'. We would like to
thank Jean-Michel Poggi (Universite Paris-Sud, France) for
providing us with the Paris electrical load consumption data.

\section*{{\sc Appendix}}

Our asymptotic results will be based on the following set
of assumptions, which we detail below before proceeding to the proofs.

\subsubsection*{Main Assumptions}

We first impose an assumption on the sample paths of the
underlying stochastic process.

\vspace*{0.3cm}

{\bf Assumption (A1):} \ The sample paths of the strictly stationary process
$Z=(Z_{i};~ i \in \NN^{+})$ are assumed to lie within a Besov space
$B^{s}_{p,q}$, where $0 < s < r$, $1 \leq p, q \leq \infty$, and $r$ is the
regularity of the scaling functions associated with the regular multiresolution
analysis discussed in Section~\ref{subsec:orthsp}.

{\bf Assumption (A2):} \ When we only observe a fixed number $P$ of samples
values from each sample path, we assume that the sampled paths of the strictly
stationary process $Z=(Z_{i};~ i \in \NN^{+})$ are continuous on $[0,\delta)$,
and that the interpolating scaling function $\phi$ of the wavelet interpolating
basis has an exponential decay.

\vspace*{0.3cm}

{\bf Assumption (A3):} \ The $\alpha_{Z}$-mixing coefficient of the strictly
stationary process $Z=(Z_{i};~ i \in \NN^{+})$ satisfies
\begin{equation} \label{eq.mix}
\sum_{m=N}^{\infty} \alpha_{Z}(m)^{1-2/l} = O(N^{-1}) \ \ \ \ \mbox{for some} \
\ \ \ l > 4.
\end{equation}
(Note that the scaling coefficients of $Z$ inherit the above mixing property,
according to the relevant discussion in Section~\ref{subsec:orthsp}.)

\vspace*{0.3cm}

We next impose some assumptions on the joint and conditional probability
density functions of the scaling coefficients $ \xi_{i}^{(J,k)}$.

\vspace*{0.3cm}

{\bf Assumption (A4):} \ $ E|\xi_{i}^{(J,k)}|^{l} < \infty$, for $ l >4$ and
every $ k=0,1,\ldots,2^{J}-1$.

%\vspace*{0.3cm}
%
%
%{\bf Assumption (A4)} \ There exist $a>0$ and $\tau>0$ such that
%$$
%\EE\left(\exp \left\{a \bigg|
%\EE\left(\xi_{i+1}^{(J,k)}\,|\,\xi_{i}^{(J,k)}\right)\bigg|^\tau\right)\right\}
%< \infty.
%$$

\vspace*{0.3cm}

{\bf Assumption (A5):} \ The joint probability density function $
f_{\xi_{i+1}^{(J,k)},\,\xi_{i}^{(J,k)}}$ of
$(\xi_{i+1}^{(J,k)},\xi_{i}^{(J,k)})$ exist, it is absolute continuous with
respect to Lebesgue measure, and it satisfies the following conditions
\begin{enumerate}
\item[($i$)] \ $ f_{\xi_{i+1}^{(j,k)},\,\xi_{i}^{(j,k)}}$ is
Lipschitz continuous, i.e.,
$$ \Big|f_{\xi_{i+1}^{(J,k)},\,\xi_{i}^{(J,k)}}(x_{1},x_{2}) -
f_{\xi_{i+1}^{(J,k)},\,\xi_{i}^{(J,k)}}(y_{1},y_{2})\Big| \leq C\,
\|(x_{1},x_{2})-(y_{1},y_{2})\|.$$

\item[($ii$)] \ The random vector of the scaling coefficients at scale $J$ admits a compactly supported
probability density function $f$ (with support $S$) which is strictly positive
and twice continuously differentiable.

\item[($iii)$] \ The conditional probability density function of $ \xi_{i+1}^{(J,k)}$
given $ \xi_{i}^{(J,k)}$ is bounded, i.e., $
 f_{\xi_{i+1}^{(J,k)}\,|\,\xi_{i}^{(J,k)}}(\cdot \mid x) \leq C < \infty$.
\end{enumerate}

\vspace*{0.3cm}

We also impose some conditions on the kernel function and the bandwidth
associated with it.

\vspace*{0.3cm}

{\bf Assumption (A6):} \  The  (univariate) kernel  $ K$ is a
bounded symmetric density on $ \RR$ satisfying $ |K(x)-K(y)| \leq
C\,|x-y|$ for all $ x,y \in \RR$. Furthermore, $ \int xK(x)dx =0$
and $ \int x^{2}K(x)dx < \infty$.

\vspace*{0.3cm}

{\bf Assumption (A7):} \ The bandwidth $ h$ satisfies $ h \rightarrow 0$ and $
nh^{2^J}\rightarrow \infty$ as $ n \rightarrow \infty$.

\vspace{0.5cm}

Let us now explain the meaning of the above assumptions. Assumptions (A1) (or
(A2)) and (A3) are quite common in times series prediction (see Bosq, 1998).
Assumptions (A4)-(A5) are essentially made on the distributional behavior of
the scaling coefficients and, therefore, are less restrictive. They are
moreover natural in nonparametric regression. Assumption (A5)-(ii) is natural
as far as it concerns the scaling coefficients since the decay of the scaling
coefficients is ensured by the approximation properties of the corresponding
transform. Moreover, it is needed for obtained uniform consistency results. The
conditions (A6)-(A7) are classical for kernel regression estimation.

\vspace{0.5cm}

\noindent {\bf Proof of Theorem~\ref{th:cons}.} \ The difference in the two
assertions of the theorem is due to the nature of the observations. In case
$(i)$, each observed segment is a time series with fixed, finite length, and we
use an interpolating wavelet transform at the appropriate resolution $J$ that
makes interpolation error negligible, i.e.,
$$
\EE({\cal P}_J(Z_{n+1})\,|\,Z_n) = \EE(Z_{n+1})\,|\,Z_n).
$$
In case $(ii)$, the observed segments are continuous-time
stochastic processes and one can not neglect anymore the
approximation error due to the projection of the observed segment
onto the scaling space $V_J$. However, under our assumptions, and
according to the results recalled in Section 2, this error is
uniformly bounded by a constant times $2^{-sJ}$, where $s$ denotes
the Besov smoothness index $s$, i.e.,
\begin{equation}\label{aux1}
||\EE({\cal P}_J(Z_{n+1})\,|\,Z_n) -  \EE(Z_{n+1})\,|\,Z_n)|| =
O\left(2^{-sJ}\right),
\end{equation}
resulting in a different rate for the second case. Hence, in both
cases, we proceed by deriving the appropriate rates for
$$
||Z_{n+1|n}^J - \EE({\cal P}_J(Z_{n+1})\,|\,Z_n)||.
$$
We first show that, as $ n \rightarrow \infty$,
\begin{equation} \label{eq.consxi}
||\Xi_{n+1|\,n}  \ - \ \EE(\Xi_{n+1}\,|\,\Xi_{n})|| \rightarrow  0, \quad
\textit{almost surely}.
\end{equation}
For this, it suffices to show that, for every $k=0,1, \ldots,
2^{J}-1$, as $ n \rightarrow \infty$,
$$
\xi_{n+1\,|\,n}^{(J,k)} \rightarrow \EE(\xi_{n+1}^{(J,k)}\,|\,\xi_{n}^{(J,k)}),
\quad \textit{almost surely}.
$$
Let $ x \in \RR^{2^{J}}$, let $ \Xi_{n+1\,|\,n}(x)$ be the value of $
\Xi_{n+1\,|\,n}$ in (\ref{eq:XXseg}) for  $ \Xi_{n}=x$, and denote by
$\xi_{n+1\,|\,n}^{(J,k)}(x)$ the $k$-th component of $ \Xi_{n+1|\,n}(x)$.
Consider the $ 2^{J}$-dimensional random variable $ W_{l}={\cal C}(\Xi_{l})$,
and denote by $ f_{\xi_{l+1}^{(J,k)}, W_{l}}$ and $ f_{W_{l}}$ the joint and
marginal densities of $ (\xi_{l+1}^{(J,k)},W_{l}) $ and $ W_{l}$, respectively.
Notice that because of (A4), and the fact that $W_{l}$ is a linear
transformation of $ \Xi_{l}$,  $ f_{\xi_{l+1}^{(J,k)}W_{l}}$ and $ f_{W_{l}}$
exist with respect to Lebesgue measure for every $k=0,1, \ldots, 2^{J}-1$. Let
$$
\widehat{f}_{W_{l}}(x)= (nh_{n}^{2^J})^{-1}\sum_{m=1}^{n-1} K(D(x,
W_{m})/h_{n})
$$
and notice that $\widehat{f}_{W_{l}}(x)$ is the kernel estimator of the
$2^{J}$-dimensional density $f_{W_{l}}(x)$. For notational convenience, in what
follows, let $\Phi_{n,k}(x)= \EE(\xi_{n+1}^{(J,k)}\,|\,\xi_{n}^{(J,k)}=x)$ and
$\hat{g}_{n,k}(x)=(nh_{n}^{2^J})^{-1}\sum_{m=1}^{n-1} K(D(x, W_{m})/h_{n})
\xi_{m+1}^{(J,k)}$.

We then have
\begin{equation} \label{eq.split1}
\left(\xi_{n+1\,|\,n}^{(J,k)}(x)
-\EE(\xi_{n+1}^{(J,k)}\,|\,\xi_{n}^{(J,k)}=x)\right) =
\frac{1}{\widehat{f}_{W_{l}}(x)}\Big\{ \hat{g}_{n,k}(x)  - \Phi_{n,k}(x)
{f}_{W_{l}}(x) \Big\} - \frac{\Phi_{n,k}(x)}{\widehat{f}_{W_{l}}(x)} \Big\{
\widehat{f}_{W_{l}}(x) -{f}_{W_{l}}(x) \Big\}.
\end{equation}
Using now the assumptions of the theorem, the above decomposition, Lemma 2.1
and Theorem 3.2 of Bosq (1998), it follows that, as $n \rightarrow \infty$,
\begin{equation}\label{arate}
\sup_{x \in S} \Big|\xi_{n+1\,|\,n}^{(J,k)}(x)
-\EE(\xi_{n+1}^{(J,k)}\,|\,\xi_{n}^{(J,k)}=x)\Big| =  O\left(\left(\frac{\log
n}{n}\right)^{1/(2+2^J)} \right), \quad \textit{almost surely}.
\end{equation}
Recalling now that our estimator is defined as
$$
Z^J_{n+1\,|\,n}(t) = \sum_{k=0}^{2^J-1} \xi_{n+1\,|\,n}^{J,k}
\phi_{J,k},
$$
using the rate given in expression~(\ref{arate}), and the fact that we have
used a regular multiresolution analysis, we have, as $n \rightarrow \infty$,
\begin{eqnarray}
\sup_{t} | Z_{n+1,n}^J(t) -  \EE(P_J(Z_{n+1}(t))\,|\,Z_n)| & = & 2^{J/2} \max_k
\bigg| \xi_{n+1|n}^{J,k} -\EE(\xi_{n+1}^{(J,k)}\,|\,\xi_{n}^{(J,k)})\bigg|
\sup_t \sum_{k=0}^{2^J-1} |\phi(2^Jt-k)| \nonumber \\
&=& O\left(2^{J/2} \left(\frac{\log n}{n}\right)^{1/(2+2^J)} \right), \quad
\textit{almost surely}. \label{eq:bfb}
\end{eqnarray}
The above bound (\ref{eq:bfb}), together with inequality~(\ref{aux1}), ensures
the validity of both the assertions $(i)$ and $(ii)$. This completes the proof
of Theorem~\ref{th:cons}. \hfill $\Box$

\vspace{0.3cm}

\noindent {\bf Proof of Theorem~\ref{th:boot}.} \ For every $ t_{i} \in
\{t_{1}, t_{2}, \ldots, t_{P}\}$, note that $ Z_{n+1}(t_i)=\xi_{n+1}^{(J,i)}$.
Since
\begin{eqnarray*}
L_{n+1,\alpha}^{\ast}(t_i) & = & R_{n+1,\alpha}^{\ast}(t_i) + Z_{n+1\,|\,n}^{J}(t_{i}) \\
& = & Z_{n+1}^{\ast}(t_i) - E(Z_{n+1}(t_{i}) \,|\, Z_{n}) +
(Z_{n+1\,|\,n}^{J}(t_{i})-E(Z_{n+1}(t_{i}) \,|\, Z_{n})),
\end{eqnarray*}
and, as $n \rightarrow \infty$, $|E(Z_{n+1}(t_{i}) \,|\,
Z_{n})-Z_{n+1\,|\,n}^{J}(t_{i})| \rightarrow 0$ \textit{in probability}, it
suffices to show that the distribution of $ Z_{n+1}^{\ast}(t_i) -
E(Z_{n+1}(t_{i}) \,|\, Z_{n})$ approximates correctly the conditional
distribution of $ Z_{n+1}(t_i) - E(Z_{n+1}(t_{i}) \,|\, Z_{n})$ given $ Z_{n}$.

Now, given $ Z_{n}=x$, i.e., given $\Xi_{n}=\widetilde{x}$, we have
\begin{eqnarray} \label{eq.kern1}
\PP(Z_{n+1}^{\ast}(t_i) - E(Z_{n+1}(t_{i}) \,|\, Z_{n}=x) \leq y)
& = & \sum_{m=1}^{n-1}1_{(-\infty,y]}(Z_{m+1}(t_{i})-E(Z_{n+1}(t_{i}) \,|\, Z_{n}=x))w_{n,m} \nonumber \\
& = & \sum_{m=1}^{n-1}1_{(-\infty,\overline{y}]}(Z_{m+1}(t_i)) w_{n,m} \nonumber \\
& = & \sum_{m=1}^{n-1}1_{(-\infty,\widetilde{y}]}(\xi_{m+1}^{(J,i)})w_{n,m} \nonumber \\
& = &  \frac{\sum_{m=1}^{n-1}{\bf
1}_{(-\infty,\widetilde{y}]}(\xi_{m+1}^{(J,i)})\, K( D({\cal C}(\widetilde{x}),
{\cal C}(\Xi_{m}))/h_{n}) }{n^{-1} + \sum_{m=1}^{n-1}K( D({\cal
C}(\widetilde{x}), {\cal
C}(\Xi_{m}))/h_{n}) } \\
& & + O(n^{-1}), \nonumber
\end{eqnarray}
where $\widetilde{y}=y+E(Z_{n+1}(t_{i}) \,|\, Z_{n}=x)$. Note that
(\ref{eq.kern1}) is a kernel estimator of the conditional mean $ E({\bf
1}_{(-\infty,\widetilde{y}]}(\xi_{n+1}^{(J,i)})\,|\,\Xi_{n}=\widetilde{x}) =
\PP(\xi_{n+1}^{(J,i)} \leq \widetilde{y}\,|\, \Xi_{n}=\widetilde{x})$, i.e., of
the conditional distribution of $ \xi_{n+1}^{(J,i)}$ given that $
\Xi_{n}=\widetilde{x} $. Denote now the conditional distribution of
$\xi_{n+1}^{(J,i)}$ given $ \Xi_{n}$ by $
F_{\xi_{n+1}^{(J,i)}\,|\,\Xi_{n}}(\cdot \,|\, \Xi_{n})$ and its kernel
estimator given in (\ref{eq.kern1}) by $
\widehat{F}_{\xi_{n+1}^{(J,i)}\,|\,\Xi_{n}}(\cdot \,|\, \Xi_{n})$. Then by the
same arguments as in Theorem~\ref{th:cons} we get that, for every $y\in \RR$,
as $n \rightarrow \infty$,
\[ \sup_{x\in S }\bigg|\widehat{F}_{\xi_{n+1}^{(J,i)}\,|\,\Xi_{n}}(y \,|\,x) -
F_{\xi_{n+1}^{(J,i)}\,|\,\Xi_{n}}( y\,|\, x)\bigg| \rightarrow 0, \quad
\textit{in probability}.\]

It remains to show that the above convergence is also uniformly over $ y$. Fix
now a $ x $ in the support $S$ of $ \Xi_{n}$, and let $ \epsilon >0$ arbitrary.
Since $ F_{\xi_{n+1}^{(J,i)}\,|\,\Xi_{n}}( y\,|\, x)$ is continuous we have
that, for every $k\in \NN$,  points  $-\infty = y_{0} < y_{1} < \ldots <
y_{k-1} < y_{k} = \infty $ exist such that $ F_{\xi_{n+1}^{(J,i)}\,|\,\Xi_{n}}(
y_{i}\,|\, x)=i/k$. For $ y_{i-1} \leq y \leq y_{i}$, and using  the
monotonicity of $ \widehat{F}_{\xi_{n+1}^{(J,i)}\,|\,\Xi_{n}}$ and $
\widehat{F}_{\xi_{n+1}^{(J,i)}\,|\,\Xi_{n}}$, we have
\begin{eqnarray*}
\widehat{F}_{\xi_{n+1}^{(J,i)}\,|\,\Xi_{n}}( y_{i-1}\,|\,
x)-F_{\xi_{n+1}^{(J,i)}\,|\,\Xi_{n}}( y_{i}\,|\, x) & \leq &
\widehat{F}_{\xi_{n+1}^{(J,i)}\,|\,\Xi_{n}}( y\,|\,
x)-F_{\xi_{n+1}^{(J,i)}\,|\,\Xi_{n}}( y\,|\, x) \\ & \leq &
\widehat{F}_{\xi_{n+1}^{(J,i)}\,|\,\Xi_{n}}( y_{i}\,|\,
x)-F_{\xi_{n+1}^{(J,i)}\,|\,\Xi_{n}}( y_{i}\,|\, x).
\end{eqnarray*}
>From this, we get
\[ \bigg|\widehat{F}_{\xi_{n+1}^{(J,i)}\,|\,\Xi_{n}}( y\,|\, x)-F_{\xi_{n+1}^{(J,i)}\,|\,\Xi_{n}}( y\,|\, x)\bigg| \leq
\sup_{i}\bigg|\widehat{F}_{\xi_{n+1}^{(J,i)}\,|\,\Xi_{n}}( y_{i}\,|\,
x)-F_{\xi_{n+1}^{(J,i)}\,|\,\Xi_{n}}( y_{i}\,|\, x)\bigg| + \frac{1}{k},\] and,
therefore,
\begin{eqnarray*}
\PP\left(\bigg|\widehat{F}_{\xi_{n+1}^{(J,i)}\,|\,\Xi_{n}}(y \,|\,x) -
F_{\xi_{n+1}^{(J,i)}\,|\,\Xi_{n}}( y\,|\, x)\bigg| > \epsilon\right) & \leq &
\PP\left(\sup_{i}\bigg|\widehat{F}_{\xi_{n+1}^{(J,i)}\,|\,\Xi_{n}}(y_{i}
\,|\,x) -
F_{\xi_{n+1}^{(J,i)}\,|\,\Xi_{n}}( y_{i}\,|\, x)\bigg| + k^{-1} > \epsilon \right)\\
& \leq &
\PP\left(\sup_{i}\sup_{x}\bigg|\widehat{F}_{\xi_{n+1}^{(J,i)}\,|\,\Xi_{n}}(y_{i}
\,|\,x) - F_{\xi_{n+1}^{(J,i)}\,|\,\Xi_{n}}( y_{i}\,|\, x)\bigg|  + k^{-1} >
\epsilon \right).
\end{eqnarray*}
Now, chose $ k$ large enough such that $ 1/k<\epsilon/2$. For such a fixed $k$,
and because, for every $y \in \RR$, as $n \rightarrow \infty$,
$$
\sup_{x}\bigg|\widehat{F}_{\xi_{n+1}^{(J,i)}\,|\,\Xi_{n}}(y \,|\,x) -
F_{\xi_{n+1}^{(J,i)}\,|\,\Xi_{n}}( y\,|\, x)\bigg| \rightarrow 0, \quad
\textit{in probability},
$$
we can choose $n$ large enough such that
\[\PP\left(\sup_{
1 \leq i\leq k}\sup_{x}\bigg|\widehat{F}_{\xi_{n+1}^{(J,i)}\,|\,\Xi_{n}}(y_{i}
\,|\,x) - F_{\xi_{n+1}^{(J,i)}\,|\,\Xi_{n}}( y_{i}\,|\, x)\bigg|   >
\epsilon/2\right) < \tau,\] for any desired $\tau$. Since $\tau$ is
independent on $y$ and $x$, the desired convergence follows. This completes the
proof of Theorem~\ref{th:boot}. \hfill $\Box$

\end{document}